\documentclass[11pt,a4paper]{article}
\usepackage[T1]{fontenc}
\usepackage{lmodern}
\usepackage{mathtools,amssymb,amsthm,mathrsfs}
\usepackage{bm,microtype,enumerate}
\usepackage{array}
\usepackage[margin=1.08in]{geometry}
\usepackage[colorlinks=true,citecolor=blue,linkcolor=blue,urlcolor=blue]{hyperref}
\newtheorem{theorem}{Theorem}[section]
\newtheorem{proposition}[theorem]{Proposition}
\newtheorem{lemma}[theorem]{Lemma}
\newtheorem{corollary}[theorem]{Corollary}
\theoremstyle{remark}
\newtheorem{remarks}[theorem]{Remarks}
\newtheorem{remark}[theorem]{Remark}
\numberwithin{equation}{section}

\newcommand{\R}{\mathbb{R}}
\newcommand{\Z}{\mathbb{Z}}
\newcommand{\Sph}{\mathbb{S}^2}
\newcommand{\cA}{\mathcal{A}}
\newcommand{\cF}{\mathcal{F}}
\newcommand{\cD}{\mathscr{D}}
\newcommand{\cH}{\mathcal{H}}
\newcommand{\cK}{\mathcal{K}}
\newcommand{\cL}{\mathcal{L}}
\newcommand{\cN}{\mathcal{N}}
\newcommand{\cO}{\mathcal{O}}
\newcommand{\cW}{\mathcal{W}}
\newcommand{\cR}{\mathcal{R}}
\newcommand{\dd}{\,\mathrm{d}}
\newcommand{\dif}{\mathrm{d}}

\newcommand{\cof}{\operatorname{cof}}
\newcommand{\Div}{\operatorname{div}}
\newcommand{\m}{\boldsymbol{m}}
\newcommand{\M}{\boldsymbol{M}}
\newcommand{\n}{\boldsymbol{n}}

\newcommand{\bs}{\boldsymbol}
\newcommand{\ein}{\boldsymbol{e}}

\newcommand{\fJ}{\mathsf{J}}
\newcommand{\fL}{\mathsf{L}}
\newcommand{\fT}{\mathsf{T}}

\newcommand{\Gt}{\widetilde{G}}
\newcommand{\p}{\partial}
\newcommand{\eps}{\varepsilon}
\DeclarePairedDelimiter{\abs}{\lvert}{\rvert} 
\DeclarePairedDelimiter{\norm}{\lVert}{\rVert} 
\DeclarePairedDelimiter{\bra}{(}{)} 
\DeclarePairedDelimiter{\pra}{[}{]} 
\DeclarePairedDelimiter{\set}{\{}{\}} 
\DeclarePairedDelimiter{\scp}{\langle}{\rangle} 
\title{Global minimality and local stability of chiral skyrmions}
\author{Xinye Li}
\AtEndDocument{%
 \par
 \medskip
 \begin{tabular}{@{}l@{}}%
School of Mathematics and Statistics, Central South University,\\ 410083 Changsha, China\\
 \textit{Email:} \texttt{xinye.li@csu.edu.cn}\\

 \end{tabular}}
\begin{document}
\maketitle
\begin{abstract}
We study isolated skyrmions in a planar chiral magnet with exchange, Dzyaloshinskii-Moriya interaction of strength $\kappa>0$, a perpendicular applied field $h>0$, and uniaxial anisotropy with coefficient $A\in\R$. For $h>|A|$ and $4\kappa\sqrt h \le 3h-A$, we prove that the minimum in the degree-$-1$ class is attained and that every minimizer is a translate of the axisymmetric skyrmion generated by a unique strictly decreasing radial profile.  The proof uses a calibration by a continuous, distributionally closed two-form and a weighted degree identity. Throughout the same parameter regime, we establish linear stability of the axisymmetric skyrmion under arbitrary perturbations. Near the translation orbit, the energy excess is comparable to the squared $H^1$ distance to that orbit.
\end{abstract}


\section{Introduction and main results}\label{sec:intro}

Magnetic skyrmions are localized magnetization textures with nonzero topological degree~\cite{Leonov2016}. In the planar chiral model considered here, the energy combines exchange and Dzyaloshinskii-Moriya interaction (DMI) with a local potential arising from a perpendicular applied field and uniaxial anisotropy. 
In two dimensions the exchange energy is invariant under dilations, leaving the size of a Belavin-Polyakov profile undetermined~\cite{BelavinPolyakov1975}. The competition between DMI and the local potential can select a characteristic skyrmion scale in near-conformal regimes~\cite{DoeringMelcher2017,GustafsonWang2021,BernandMantel2021}.
An axisymmetric configuration is a natural candidate for an isolated skyrmion. We ask when this candidate minimizes the energy among all configurations of the same degree and whether the energy excess measures the size of a deformation after the center has been allowed to move. 

We prove existence and classify all degree-$-1$ global minimizers for a two-dimensional chiral magnet with Zeeman coupling and uniaxial anisotropy. In an explicit parameter regime, their smooth representatives are exactly the translates of the axisymmetric map generated by a unique center-regular, strictly decreasing radial profile. We also establish a two-sided quadratic bound of the energy excess in terms of the $H^1$ distance to the translation orbit, valid in a sufficiently small $H^1$ neighborhood of that orbit.

We consider magnetization fields $\m=(m',m_3):\R^2 \to \Sph$, where $m'=(m_1,m_2)$, with the background state $\ein_3=(0,0,1)$. Using the interfacial DMI convention, for $\kappa, h>0$ and $A \in \R$, we define 
\begin{equation}\label{eq:energy}
E(\m)=\int_{\R^2} \frac{1}{2}|D \m|^2 \dd x + \kappa \cD(\m) + \int_{\R^2} W(\m) \dd x
\end{equation}
where
\[
\cD(\m)=-2 \int_{\R^2} m' \cdot \nabla m_3 \dd x \quad\text{and}\quad
W(\m)=h(1-m_3) + \frac{A}{2}(1-m_3^2).
\]
Positive $A$ denotes easy-axis anisotropy and negative $A$ denotes easy-plane anisotropy. The natural admissible space is
\[
\cA=\set*{\m:\R^2 \to \Sph \text{ measurable with } \m-\ein_3 \in H^1(\R^2;\R^3)},
\]
which ensures all terms in the energy~\eqref{eq:energy} are finite. 
In this work, we study the minimizers of $E$ on the class $\cA_{-1} = \set{\m \in \cA \text{ with } Q(\m)={-1}}$, where $Q$ denotes the topological degree
\[
Q(\m)=\frac{1}{4\pi} \int_{\R^2} \m \cdot \p_1 \m \times \p_2 \m \dd x.
\]
The degree $Q$ is integer-valued and continuous under strong affine $H^1$ convergence on $\cA$; see~\cite[\S~2(d)--(e)]{Melcher2014}.

Let $R_\gamma \in SO(2)$ denote planar rotation by $\gamma$, and let $\bs R_\gamma=\mathrm{diag}(R_\gamma,1)$ denote the target rotation about the $\ein_3$-axis. For $z\in \R^2 \setminus \set{0}$, set
\[
\ein_r(z) = \frac{z}{|z|} \quad\text{and}\quad \ein_{\varphi}(z) = R_{\frac{\pi}{2}} \ein_r(z).
\]
We now state the main result.
\begin{theorem}\label{thm:main}
Fix $\kappa,h>0$ and $A \in \R$ satisfying
\begin{equation}\label{eq:range}
h>|A|,\qquad 4\kappa\sqrt h \le 3h-A.
\end{equation}
There exists a unique profile $\theta \in C([0,\infty)) \cap C^\infty((0,\infty))$ satisfying
\begin{equation}\label{eq:radial-ode}
(r\theta')'=\frac{\sin \theta \cos \theta}{r} - 2\kappa \sin^2 \theta + hr \sin \theta + Ar \sin \theta \cos \theta,
\end{equation}
with
\begin{equation}\label{eq:radial-bc}
\theta(0)=\pi, \quad \theta(\infty)=0, \quad 0<\theta(r)<\pi, \quad \theta'(r)<0 \quad\text{for all } r>0,
\end{equation}
and such that the associated map $\M$ defined below extends smoothly across the origin.
Define 
\[
\M(x)=
\begin{cases}
\bra*{\sin \theta(|x|) \ein_r(x), \cos \theta(|x|)}, & x \neq 0,\\
-\ein_3, &x=0
\end{cases}
\]
and $\M_c(x)=\M(x-c)$ for $c\in\R^2$. Then $\M_c\in\cA_{-1}$ is smooth and globally minimizes $E$ over $\cA_{-1}$. Conversely, every global minimizer of $E$ over $\cA_{-1}$ equals $\M_c$ almost everywhere for a unique $c\in\R^2$.
\end{theorem}

\begin{corollary}\label{coro:coercivity}
Assume the hypotheses of Theorem~\ref{thm:main}.
\begin{enumerate}[(i)]
\item The centered minimizer $\M$ is linearly stable modulo translations.
More precisely, its Hessian $\cH_{\M}$ satisfies
\[
\cH_{\M}(\bs\eta,\bs\eta)\ge \lambda\norm{\bs\eta-\bs\eta_0}_{H^1}^2
\]
for some $\lambda>0$ and every $\bs\eta\in H^1(\R^2;T_{\M}\Sph)$.
Here $\bs\eta_0$ is the $L^2$-orthogonal projection of $\bs\eta$ onto $\mathrm{span}\{\p_1\M,\p_2\M\}$, the kernel of $\cH_{\M}$ on the tangent space.

\item The minimizer is locally energetically stable modulo translations.
Let $\cO=\set*{\M_c:\,c\in\R^2}$ and define
\[
d_{H^1}(\m,\cO):=\inf_{c\in\R^2}\norm{\m-\M_c}_{H^1}.
\]
For fixed parameters, there exist $\eps_*,c_1,c_2>0$ such that every $\m\in\cA_{-1}$ with $d_{H^1}(\m,\cO)<\eps_*$ satisfies
\begin{equation}\label{eq:local-coercivity}
c_1 d_{H^1}(\m,\cO)^2 \le E(\m)-E(\M) \le c_2 d_{H^1}(\m,\cO)^2.
\end{equation}
\end{enumerate}
\end{corollary}

\bigskip

\begin{remarks}\label{rem:intro}
\begin{enumerate}[(1)]
\item Assume $\kappa>0$ and $h>|A|$, and let $\theta$ be a regular solution of~\eqref{eq:radial-ode} satisfying~\eqref{eq:radial-bc}. The density constructed in Section~\ref{sec:construction} gives a pointwise lower bound for the energy density if and only if $(-\theta'(0))^2\ge h$. 
At equality the quadratic remainder degenerates at the core, and the associated map still minimizes $E$ in $\cA_{-1}$, uniquely up to translations. 
The remainder is uniformly positive definite precisely when the inequality is strict, as guaranteed by~\eqref{eq:range}.
\item The uniqueness of $\theta$ in Theorem~\ref{thm:main} is asserted only within the class of center-regular solutions satisfying~\eqref{eq:radial-bc}. The theorem does not exclude other regular nonmonotone solutions of~\eqref{eq:radial-ode} with the same endpoint values. However, if such solutions exist and have finite energy, their associated maps cannot be additional global minimizers in $\cA_{-1}$.
\item The interfacial DMI convention aligns the Néel profile with the radial frame. In the bulk convention with the same normalization,
\[
\cD_{\mathrm{bulk}}(\m)
=2\int_{\R^2}\bra*{m_1\p_2m_3-m_2\p_1m_3}\dd x
=\cD(\bs R_{-\pi/2}\m).
\]
Target rotations about $\ein_3$ preserve the other energy terms, the degree, and $H^1$ distances. Thus Theorem~\ref{thm:main} and Corollary~\ref{coro:coercivity} carry over with $\bs R_{\pi/2}\M$ as the bulk reference map. In either convention, rotation by $\pi$ reverses the DMI sign. The case $\kappa<0$ therefore follows from the positive case with $|\kappa|$ under $h>|A|$ and $4|\kappa|\sqrt h\le 3h-A$.
\end{enumerate}
\end{remarks}
\paragraph{Relation to previous rigorous work}
Existence of minimizers of degree-$-1$ on the whole plane has been established in several parameter regimes. In the pure Zeeman model, the infimum in the degree-$-1$ class is attained for $h>\kappa^2$~\cite{Melcher2014}. 
Existence also holds for sufficiently large fields when $A>0$~\cite{HoffmannEtAl2017}, and in the zero-field easy-axis model when $A \ge 8\kappa^2$ in our normalization~\cite{BernandMantel2021}. Conformal-limit analyses describe the approach to the Belavin--Polyakov family and the selection of scale and helicity \cite{DoeringMelcher2017,GustafsonWang2021,BernandMantel2020,BernandMantel2021}.
A quantitative rigidity estimate controls the distance of degree-$\pm1$ maps to that family~\cite{BernandMantel2021}. Asymptotic matching has also been used to study the profiles and radii of axisymmetric skyrmions, including when an applied field is combined with easy-axis anisotropy~\cite{KomineasMelcherVenakides2023}.

For the square-potential family $A=-h$, a Bogomol'nyi decomposition at $h=\kappa^2$ yields explicit solutions~\cite{BartonSinger2020}.
Stability thresholds and fixed-degree minimization in this family have been studied in~\cite{IbrahimShimizu2023,IbrahimMiuraRomanShimizu2025}.
For bounded films with constant Dirichlet data, existence and concentration of skyrmion minimizers have been established in the weak DMI regime~\cite{MonteilMuratovSimonSlastikov2023}. 
Higher-degree minimizers exist on sufficiently large or slender domains, and their concentration has been analyzed in a strong-anisotropy regime~\cite{MuratovSimonSlastikov2025}.

In the zero-field easy-axis model ($h=0$, $A>0$), Di Fratta et al.~\cite{DiFrattaMuratovSimonSlastikov2026} proved that, for $0<\kappa<2\sqrt{A}/\pi$ in our normalization, the minimization problem on $\mathbb{R}^2$ in the degree-$-1$ class admits an axisymmetric minimizer, unique up to translations.
At $\kappa=2\sqrt{A}/\pi$, the infimum is zero and is not attained, and for larger $\kappa$ the energy is unbounded below.
For further symmetry and uniqueness results in related variational problems, see the references cited there.

At sufficiently large fields in the pure Zeeman model, a spectral gap for the Hessian yields a quadratic $H^1$ lower bound on the energy excess near the translation orbit of the axisymmetric skyrmion~\cite{LiMelcher2018}. 
For the combined Zeeman--anisotropy model, exponentially decaying equivariant solutions exist whenever $h>|A|$, and the Hessian at the constructed solutions is nonnegative on all $H^1$ tangent fields if $h-A \ge 8\kappa^2$ in our normalization~\cite{IbrahimShimizu2026}. 

For the combined Zeeman--anisotropy potential, we establish two global singular-value barriers and a strict center-slope bound under~\eqref{eq:range}. 
Together with the pole limits, these estimates make the quadratic calibration remainder uniformly positive definite on the whole target sphere, including the skyrmion core, for each fixed parameter triple in this regime. 
The corresponding remainder in the zero-field easy-axis calibration degenerates at the core; see~\cite[Proposition~6.1]{DiFrattaMuratovSimonSlastikov2026}.
This uniform positivity is the key to global rigidity and quantitative stability modulo translations.
\paragraph{The calibration argument}
We seek a pointwise lower bound for the energy density that is attained by the radial reference map and has the same integral for all maps of a given degree. Before the divergence correction, the energy density can be written as
\[
f(\n,B)=\frac12|B|^2+\kappa G(\n):B+W(\n)
\]
where $G$ is the interfacial DMI coefficient. For $\n\ne\ein_3$, let $\fJ(\n)$ be the gradient of the radial reference map at the point where its magnetization equals $\n$. Its extension to $\ein_3$ is established in Section~\ref{sec:construction}. Sharpness at $\fJ$ requires matching $f$ in value and first $B$-derivative. We first seek a linear term $p(\n):B$ with zero integral.
A divergence with suitable decay has this property and is linear in $D\m$ by the chain rule $\Div(Y\circ\m)=(\dif Y(\m))^{\fT}:D\m$ for $Y:\Sph\to\R^2$.
With this term alone, derivative matching fixes $p$ and hence $p:\fJ$.

A correction depending only on $\n$ could match the value, but its integral scales with the area under a degree-preserving dilation. The degree counts target preimages with sign, which suggests using the Jacobian, the density of oriented target area. 
Weighting it by a continuous function on the target gives an adjustable coefficient for value matching and an integral that still depends only on the degree. We therefore take the candidate calibrating density to be
\[
\ell(\n,B)=p(\n):B+\beta(\n)\det B.
\]
With the linear term realized as a divergence, the calibrating density has a weighted-Jacobian and divergence structure also present in the zero-field easy-axis calibration of~\cite{DiFrattaMuratovSimonSlastikov2026}.

Matching the first $B$-derivative and the value at $\fJ$ gives
\begin{align*}
&p=D_Bf(\n,\fJ)-\beta\cof\fJ,\\
&\beta\det\fJ=D_Bf(\n,\fJ):\fJ-f(\n,\fJ)=\frac12|\fJ|^2-W(\n).
\end{align*}
The radial equation gives the compatibility condition for the divergence potential; see Remark~\ref{rem:stress-closure} for its interpretation through stress conservation. For the whole-plane comparison, we subtract $\kappa\Div m'$ from both densities.
Lemma~\ref{lem:closedness} provides the pole extensions and decay estimates ensuring that the corrected divergence is absolutely integrable and has zero integral.
We express this structure as the graph pullback of a continuous, distributionally closed two-form on $\R^2\times\Sph$.
Its target part is the weighted sphere-area form and its mixed part is exact. 
This representation connects the construction with first-order null Lagrangians~\cite{OlverSivaloganathan1988} and energy calibrations on compact  manifolds~\cite{Hattori2024}.

The common DMI correction cancels in $f-\ell$. With $C=B-\fJ$, the contact conditions and $2|\det C|\le |C|^2$ give
\[
f(\n,B)-\ell(\n,B)=\frac12|C|^2-\beta\det C
\ge \frac12\bigl(1-\norm{\beta}_{L^\infty(\Sph)}\bigr)|C|^2.
\]
For each parameter triple satisfying~\eqref{eq:range}, the two singular-value barriers and the pole limits give $|\beta|<1$ on the entire sphere, including the skyrmion core. By continuity and compactness, $\norm{\beta}_{L^\infty(\Sph)}<1$.
This uniform gap, together with the weighted degree identity, bounds the energy excess $E(\m)-E(\M)$ from below by a positive multiple of $\norm{D\m-\fJ(\m)}_{L^2}^2$ throughout $\cA_{-1}$.
In Section~\ref{sec:coercivity}, we linearize this residual to identify the translation kernel and then use a far-field estimate and local compactness to obtain the spectral gap.
\paragraph{Outline of the paper}
Section~\ref{sec:pre} gives the integrable energy-space formulation and the first- and second-variation identities. Section~\ref{sec:criterion} proves the calibration criterion. Sections~\ref{sec:radial} and~\ref{sec:barriers} establish the monotone radial profile and the two barriers. Section~\ref{sec:construction} classifies the global minimizers. Section~\ref{sec:coercivity} proves linear stability and derives a local energy estimate modulo translations.

Throughout the paper, primes on scalar functions of $r$ denote differentiation with respect to $r$. For a vector field $\m=(m_1,m_2,m_3)$, we write $m'=(m_1,m_2)$ for its in-plane component. Derivatives with respect to other variables are indicated by subscripts or explicit differential operators.

\section{Preliminaries}\label{sec:pre}
\subsection{Intrinsic determinants and cofactors}
Let $e_1,e_2$ denote the standard basis of $\R^2$, and let $\sigma$ be the standard area-form on $\Sph$ normalized by $\int_{\Sph} \sigma = 4\pi.$ For $\n\in\Sph$ and $B:\R^2\to T_{\n}\Sph$, we use the intrinsic determinant
\[
 \det B=\sigma_{\n}(Be_1,Be_2).
\]
After choosing a positively oriented orthonormal frame of $T_{\n}\Sph$, the map $B$ is represented by a $2\times 2$ matrix. We define the intrinsic cofactor $\cof B:\R^2\to T_{\n}\Sph$ by
\[
\frac{\dif}{\dif t} \Big|_{t=0} \det(B+tC)=\cof B:C \quad\text{for every } C:\R^2 \to T_{\n}\Sph,
\]
where $:$ denotes the inner product induced by the Euclidean metrics on $\R^2$ and $T_{\n}\Sph$. In any positively oriented orthonormal frame of $T_{\n}\Sph$, the map $\cof B$ is represented by the usual cofactor matrix of the matrix representing $B$. Thus $\cof B$ is independent of the chosen frame, although its matrix representation depends on that choice. In particular,
\begin{equation}\label{eq:2D-algebra}
 |\cof B|=|B|, \qquad |\det B|\le\frac{1}{2}|B|^2.
\end{equation}
\subsection{A null-Lagrangian adjustment of DMI}
The DMI density $-2m'\cdot\nabla m_3$ is absolutely integrable for every $\m\in\cA$. For the pointwise energy comparison, we now represent the same functional using a tangent-valued coefficient paired with $D\m$.

For $\n=(n',n_3)\in\Sph$, define the linear map $G(\n):\R^2\to T_{\n}\Sph$ by
\[
G(\n)e_i=n_3\ein_i-n_i\ein_3, \qquad i=1,2,
\]
and let $\pi(\n)=n'$ be the in-plane projection. The usual interfacial Lifshitz density is
\[
G(\m):D\m =m_3\Div m'-m'\cdot\nabla m_3
=\Div m'+(m_3-1)\Div m'-m'\cdot\nabla m_3.
\]
The linear divergence $\Div m'$ need not belong to $L^1(\R^2)$ under the sole assumption $\m-\ein_3\in H^1(\R^2;\R^3)$. Since
\[
(\dif\pi(\n))^{\fT}e_i=\ein_i-n_i\n,
\qquad
(\dif\pi(\m))^{\fT}:D\m=\Div m',
\]
where $\fT$ denotes the Euclidean adjoint, we introduce the modified coefficient
\[
\Gt=G-(\dif\pi)^{\fT},
\qquad \Gt(\n)e_i=(n_3-1)\ein_i+n_i(\n-\ein_3).
\]
We verify below that this coefficient yields an absolutely integrable density whose integral equals the original DMI functional $\cD(\m)$.

The bound $|\Gt(\n)|\lesssim|\n-\ein_3|$ implies $\Gt(\m)\in L^2$ and $\Gt(\m):D\m\in L^1$ for every $\m\in\cA$. Moreover,
\[
\Gt(\m):D\m+2m'\cdot\nabla m_3=\Div\big[(m_3-1)m'\big]
\]
in the sense of distributions.
For any $Z\in W^{1,1}(\R^2;\R^2)$, choose cutoff functions $\chi_R \in C_c^\infty(B_{2R})$  satisfying
\[
0 \le \chi_R \le 1, \qquad \chi_R \equiv 1 \text{ on } B_R, \qquad \abs*{D\chi_R} \le \frac{C}{R},
\]
where $C>0$ is a constant independent of $R$. Integration by parts gives
\[
\abs*{\int_{\R^2} \chi_R \Div Z \dd x} = \abs*{\int_{\R^2} Z \cdot D\chi_R \dd x } 
\le \frac{C}{R} \int_{B_{2R} \setminus B_R} |Z| \to 0 \quad\text{as } R \to \infty.
\]
Since $\Div Z\in L^1$, dominated convergence yields
\begin{equation}\label{eq:zero-divergence}
 \int_{\R^2}\Div Z\,\dd x=0.
\end{equation}
For $Z=(m_3-1)m'$, the product rule gives
\[
|Z|\le|\m-\ein_3|^2,\qquad
|DZ|\le2|\m-\ein_3|\,|D\m| \quad\text{a.e.}
\]
Since $\m-\ein_3\in H^1(\R^2;\R^3)$, both right-hand sides belong to $L^1(\R^2)$ by H\"older's inequality.
Thus $Z\in W^{1,1}(\R^2;\R^2)$.
Applying~\eqref{eq:zero-divergence} to $Z$ and using the preceding divergence identity yields
\[
\int_{\R^2}\Gt(\m):D\m\,\dd x =-2\int_{\R^2}m'\cdot\nabla m_3\,\dd x =\cD(\m).
\]
Consequently, the original energy~\eqref{eq:energy} has the equivalent representation
\begin{equation}\label{eq:gauged-energy}
E(\m)=\int_{\R^2}
\frac12|D\m|^2+\kappa\Gt(\m):D\m+W(\m)\dd x.
\end{equation}
Every term in this representation is absolutely integrable, and this is the form used in the calibration argument.
\subsection{First and second variations}\label{sec:secondVariation}
As in \cite{LiMelcher2018}, we express the energy and its second variation using a bounded symmetric form on $H^1(\R^2;\R^3)$. For $\bs\eta, \bs\zeta \in H^1(\R^2;\R^3)$, written as $\bs\eta=(\eta',\eta_3)$ and $\bs\zeta=(\zeta',\zeta_3)$, define the bounded symmetric bilinear form
\[
\cH_\infty(\bs\eta,\bs\zeta)=\int_{\R^2} D\bs\eta:D\bs\zeta-2\kappa\bra*{\eta'\cdot \nabla\zeta_3 + \zeta' \cdot \nabla \eta_3} +(h+A)\eta' \cdot \zeta' + h\eta_3 \zeta_3 \dd x.
\]
We use the notation $\cH_\infty(\bs\eta)=\cH_\infty(\bs\eta,\bs\eta)$ and the same convention for $\cH_{\M}$ below. For $\m \in \cA$, the identities
\[
1-m_3 = \frac{1}{2}|\m-\ein_3|^2 \quad\text{and}\quad 1-m_3^2 = |m'|^2
\]
give
\begin{equation}\label{eq:EHessian}
E(\m)=\frac{1}{2}\cH_\infty(\m-\ein_3).
\end{equation}
Let $\M \in \cA$. For $\bs\eta\in  L^\infty \cap H^1(\R^2; T_{\M}\Sph)$, consider the variation
\[
\m_t = \frac{\M +t\bs\eta}{|\M +t\bs\eta|}.
\]
We call $\M$ a critical point of $E$ if the first variation vanishes for every such $\bs\eta$.
By~\eqref{eq:EHessian}, this means
\[
\frac{\dif}{\dif t}E(\m_t)\Big|_{t=0} = \cH_\infty(\M-\ein_3,\bs\eta)=0.
\]
Assume from now on that $\M$ is a critical point of $E$ and that $\M-\ein_3\in H^2(\R^2;\R^3)$. Integration by parts gives
\[
\cH_\infty(\M-\ein_3,\bs\eta)=\int_{\R^2} \bs F \cdot \bs\eta \dd x,
\quad\text{where}\quad
\bs F=\begin{pmatrix}
-\Delta M' - 2\kappa \nabla M_3 + (h+A)M'\\
-\Delta M_3 + 2\kappa \Div M' + h(M_3-1)
\end{pmatrix}.
\]
We denote the pointwise orthogonal projection onto $T_{\M}\Sph$ by $P_{\M} = I - \M \otimes \M$. For $\bs\zeta\in C_c^\infty(\R^2;\R^3)$, the field $P_{\M}\bs\zeta$ is a bounded $H^1$ tangent field. Testing criticality with it gives $P_{\M} \bs F=0$. Consequently,
\begin{equation}\label{eq:Hinfity}
\cH_\infty(\M-\ein_3,\bs\eta) = \int_{\R^2} \Lambda(\M) \M \cdot \bs \eta \dd x \quad\text{for every}\quad
\bs\eta \in H^1(\R^2;\R^3),
\end{equation}
where
\[
\Lambda(\M) = \M \cdot \bs F = |D\M|^2+2\kappa (M_3 \Div M' - M' \cdot \nabla M_3) + h(1-M_3)+A|M'|^2.
\]
In particular, $D\M \in L^4$ by the Gagliardo-Nirenberg inequality, and thus $\Lambda(\M) \in L^2$.

Set
\[
\m_{s,t}=\frac{\M+s\bs\eta+t\bs\zeta}{|\M+s\bs\eta+t\bs\zeta|} \quad\text{for}\quad \bs\eta,\bs\zeta\in H^2(\R^2;T_{\M}\Sph).
\]
Applying the chain rule to~\eqref{eq:EHessian} gives
\[
\frac{\p^2}{\p s\p t}\Big|_{s,t=0} E(\m_{s,t}) =\cH_\infty(\bs \eta, \bs \zeta) - \cH_\infty(\M-\ein_3, (\bs \eta \cdot\bs \zeta)\M).
\]
By~\eqref{eq:Hinfity} the Hessian of $E$ at $\M$ is therefore
\begin{equation}\label{eq:Hessian}
\cH_{\M}(\bs \eta, \bs \zeta) = \cH_\infty(\bs \eta, \bs \zeta) - \int_{\R^2} \Lambda(\M) \bs\eta \cdot \bs\zeta \dd x.
\end{equation}
H\"older's inequality and $H^1(\R^2) \hookrightarrow L^4(\R^2)$ give
\[
\abs*{\int_{\R^2} \Lambda(\M) \bs\eta \cdot \bs\zeta  \dd x} 
\le \norm{\Lambda(\M)}_{L^2}\norm{\bs\eta}_{L^4}\norm{\bs\zeta}_{L^4}
\le \norm{\Lambda(\M)}_{L^2}\norm{\bs\eta}_{H^1}\norm{\bs\zeta}_{H^1}.
\]
Hence~\eqref{eq:Hessian} defines a bounded bilinear form on $H^1(\R^2;\R^3)$ and can be evaluated on the difference $\bs\xi=\m-\M$ for any $\m\in\cA$.
Expanding~\eqref{eq:EHessian} and using $\M\cdot\bs\xi=-\frac{1}{2}|\bs\xi|^2$ and~\eqref{eq:Hinfity} yields
\begin{equation}\label{eq:Edifference}
E(\m)-E(\M)=\frac{1}{2}\cH_\infty(\bs\xi) + \cH_\infty(\M-\ein_3,\bs\xi)=\frac{1}{2}\cH_{\M}(\bs\xi).
\end{equation}

\section{A calibration criterion for global minimality and rigidity}\label{sec:criterion}
We first show that the integral of the calibrating density depends only on
the degree. We then use a pointwise defect identity to prove global
minimality and characterize equality.

\begin{lemma}\label{lem:degree}
 Let $\m \in \cA$ and $\beta \in C(\Sph)$. With
 $\m^*\sigma=(\det D\m)\,\dif x_1\wedge\dif x_2$, we have
 \begin{equation}\label{eq:degree}
 \int_{\R^2} \beta(\m) \m^* \sigma = Q(\m) \int_{\Sph} \beta \sigma.
 \end{equation}
\end{lemma}
\begin{proof}
By \cite[\S~2(d)--(e)]{Melcher2014}, there are smooth sphere-valued maps
$\m_k$, equal to $\ein_3$ outside a compact set, such that
\begin{equation}\label{eq:sphere-density}
\m_k-\ein_3\longrightarrow \m-\ein_3 \qquad\text{strongly in }H^1(\R^2;\R^3).
\end{equation}
Each $\m_k$ extends smoothly to the one-point compactification. Applying
the classical degree formula (see also \cite[Property~4]{BrezisNirenberg1995})
and then approximating $\beta$ uniformly by smooth functions gives
\begin{equation}\label{eq:smooth-weighted-degree}
 \int_{\R^2}\beta(\m_k)\,\m_k^*\sigma =\deg(\m_k)\int_{\Sph}\beta\sigma,
 \qquad \deg(\m_k)=Q(\m_k).
\end{equation}
Set $F(\n)=\beta(\n)\n$, $\boldsymbol a_k=\p_1\m_k\times\p_2\m_k$ and
$\boldsymbol a=\p_1\m\times\p_2\m$. By H\"older's inequality and~\eqref{eq:sphere-density},
$\boldsymbol a_k\to\boldsymbol a$ strongly in $L^1$. Along a subsequence
$\m_k\to\m$ almost everywhere, so boundedness and continuity of $F$
give
\[
 \|F(\m_k)\cdot\boldsymbol a_k-F(\m)\cdot\boldsymbol a\|_{L^1}
 \le \|F\|_{L^\infty}\|\boldsymbol a_k-\boldsymbol a\|_{L^1}
  + \int_{\R^2}|F(\m_k)-F(\m)|\,|\boldsymbol a|\dd x
 \longrightarrow 0,
\]
where the last integral tends to zero by dominated convergence. By
continuity of the degree under strong affine $H^1$ convergence,
$Q(\m_k)\to Q(\m)$. Passing to the limit in~\eqref{eq:smooth-weighted-degree}
proves~\eqref{eq:degree}.
\end{proof}
The next lemma combines the weighted Jacobian with a divergence term
that is linear in $D\m$.
\begin{lemma}\label{lem:graph-integral}
 Let $\beta \in C(\Sph)$, $Y \in C^1(\Sph,\R^2)$, and suppose that
 \begin{equation}\label{eq:Yt}
 |Y(\n)| \le C|\n-\ein_3|^2, \qquad |\dif Y(\n)| \le C|\n-\ein_3|.
 \end{equation}
Define
 \[
\Omega= \dif Y_1 \wedge \dif x_2 - \dif Y_2 \wedge \dif x_1 + \beta \sigma \quad\text{on }\R_x^2\times\Sph.
\] 
Then $\Omega$ is a continuous, distributionally closed two-form. For $\m\in\cA$, let
\[
 \Gamma_{\m}:\R^2\to\R^2\times\Sph,
 \qquad \Gamma_{\m}(x)=(x,\m(x))
\]
denote the graph map. Define its pullback of $\Omega$ almost everywhere by
\begin{equation}\label{eq:graph-pullback}
 \Gamma_{\m}^*\Omega
 =\pra*{ (\dif Y(\m))^{\fT}:D\m  +\beta(\m)\det D\m} \,\dif x_1\wedge \dif x_2.
\end{equation}
Both terms in~\eqref{eq:graph-pullback} are absolutely integrable, and
\begin{equation}\label{eq:graph-integral}
\int_{\R^2} \Gamma_{\m}^* \Omega= Q(\m)\int_{\Sph}\beta \sigma.
\end{equation}
\end{lemma}
\begin{proof}
The mixed part is exact, since it equals $\dif(Y_1\,\dif x_2-Y_2\,\dif x_1)$.
The term $\beta\sigma$ is pulled back from a top-degree form on $\Sph$
and is therefore distributionally closed.
Evaluating $\Omega$ on the weak derivatives of the graph map
\[
D\Gamma_{\m}(e_i)=(e_i,\p_i\m),\qquad i=1,2
\]
gives~\eqref{eq:graph-pullback}. For $Z=Y\circ\m$, the Sobolev chain rule
identifies the mixed density with $\Div Z=(\dif Y(\m))^{\fT}:D\m$.
By~\eqref{eq:Yt},
\[
|Z| \le C|\m-\ein_3|^2 \in L^1, \qquad |DZ| \le C \abs{\m-\ein_3} \, \abs*{D\m} \in L^1.
\]
Thus $Z\in W^{1,1}(\R^2;\R^2)$, so the mixed term is integrable and its
integral vanishes by~\eqref{eq:zero-divergence}. The weighted Jacobian is
integrable by~\eqref{eq:2D-algebra}, and~\eqref{eq:degree} gives its
integral, proving~\eqref{eq:graph-integral}.
\end{proof}
Matching the densities at the reference gradient leaves a quadratic
remainder that we can integrate using the fixed-degree identity.
\begin{proposition}\label{prop:calibration}
Let $\cW\in C(\Sph)$. For each $\n\in\Sph$, let $\fL(\n):\R^2\to T_{\n}\Sph$ be a linear map.
Define
\[
\begin{aligned}
& f(\n,B)=\frac12|B|^2+\fL(\n):B+\cW(\n) \quad\text{for}\quad \n\in\Sph, \, B:\R^2\to T_{\n}\Sph,\\
\qquad
& \cF(\m)=\int_{\R^2} f(\m,D\m)\,\dd x \quad\text{for}\quad \m\in\cA_{-1}.
\end{aligned}
\]
Let $\M:\R^2\to\Sph\setminus\{\ein_3\}$ be a $C^1$ diffeomorphism belonging to $\cA_{-1}$, with inverse $X$. Assume the following.
\begin{enumerate}[(i)]
\item The reference gradient $\fJ(\n)=D\M(X(\n))$ extends continuously to $\Sph$ and satisfies $|\fJ(\n)|\le C|\n-\ein_3|$.
\item There exist $\beta\in C(\Sph)$ and $Y\in C^1(\Sph;\R^2)$ such that the density
\[
\ell(\n,B)=(\dif Y(\n))^{\fT}:B+\beta(\n)\det B
\]
agrees with $f$ in value and in its first $B$-derivative at the reference gradient, i.e. for every $\n\in\Sph$,
\[
f(\n,\fJ(\n))=\ell(\n,\fJ(\n)),
\qquad
D_B f(\n,\fJ(\n))=D_B\ell(\n,\fJ(\n)).
\]
In addition,
\begin{equation}\label{eq:criterion-data1}
\gamma:=1-\norm{\beta}_{L^\infty(\Sph)}>0,\qquad
|Y(\n)|\le C|\n-\ein_3|^2,\qquad
|\dif Y(\n)|\le C|\n-\ein_3|.
\end{equation}
\end{enumerate}
Then $\cF(\m)$ is absolutely convergent for every $\m\in\cA_{-1}$.
Writing $C_{\m}=D\m-\fJ(\m)$, we have
\begin{equation}\label{eq:criterion-defect}
\cF(\m)-\cF(\M)
=\int_{\R^2}\pra*{\frac12|C_{\m}|^2-\beta(\m)\det C_{\m}} \dd x
\ge\frac{\gamma}{2}\|C_{\m}\|_{L^2}^2.
\end{equation}
Moreover, $\cF(\m)=\cF(\M)$ if and only if $\m(x)=\M(x-c)$
almost everywhere for some $c\in\R^2$.
\end{proposition}
\begin{proof}
The derivative condition determines $(\dif Y)^{\fT}$ because the $B$-derivative of the determinant at $\fJ$ is $\cof\fJ$. Substituting this expression into the value condition and using $\cof\fJ:\fJ=2\det\fJ$ determines $\cW$. 
The contact conditions are thus equivalent to
\begin{equation}\label{eq:criterion-data}
(\dif Y)^{\fT}=\fJ+\fL-\beta\cof\fJ,
\qquad
\cW=\frac12|\fJ|^2-\beta\det\fJ.
\end{equation}

These identities first give the bounds needed to define $\cF$.
By (i)-(ii) and~\eqref{eq:2D-algebra},
\begin{align*}
|\cW(\n)| &\le C|\fJ(\n)|^2\le C|\n-\ein_3|^2,\\
|\fL(\n)| &\le |\dif Y(\n)|+
(1+|\beta(\n)|)|\fJ(\n)|\le C|\n-\ein_3|.
\end{align*}
Since $\m-\ein_3\in H^1$, it follows that $\cW(\m)\in L^1$ and $\fL(\m)\in L^2$. By Cauchy--Schwarz, $\fL(\m):D\m\in L^1$, so every term in $f(\m,D\m)$ is absolutely integrable.

Let $\Omega$ be the two-form associated with $\beta$ and $Y$ in Lemma~\ref{lem:graph-integral}. Its graph pullback has density $\ell(\m,D\m)$, and for every $\m\in\cA_{-1}$,
\begin{equation}\label{eq:criterion-common-graph-integral}
\int_{\R^2}\ell(\m,D\m)\,\dd x
=\int_{\R^2}\Gamma_{\m}^*\Omega
=-\int_{\Sph}\beta\sigma.
\end{equation}

For the pointwise comparison, the identities in~\eqref{eq:criterion-data} cancel the constant and linear terms in $B-\fJ(\n)$. Using
\[
\det(B-\fJ)=\det B-\cof\fJ:B+\det\fJ,
\]
we obtain
\begin{equation}\label{eq:pointwise}
f(\n,B)-\ell(\n,B)=\frac12|B-\fJ(\n)|^2
-\beta(\n)\det(B-\fJ(\n)).
\end{equation}
For $C=B-\fJ(\n)$, the determinant bound~\eqref{eq:2D-algebra} gives
\[
\frac12|C|^2-\beta(\n)\det C
\ge\frac{1-|\beta(\n)|}{2}|C|^2
\ge\frac{\gamma}{2}|C|^2.
\]
For $\m\in\cA_{-1}$, $C_{\m}\in L^2$ by (i). Since $D\M=\fJ(\M)$, we have $f(\M,D\M)=\ell(\M,D\M)$. 
The two calibration integrals agree by~\eqref{eq:criterion-common-graph-integral}, so integrating the pointwise identity and bound gives~\eqref{eq:criterion-defect}.

Suppose that $\cF(\m)=\cF(\M)$. The defect bound gives
\begin{equation}\label{eq:equality-system}
D\m=\fJ(\m)\qquad\text{almost everywhere}.
\end{equation}
Since $\fJ$ is bounded on $\Sph$, \eqref{eq:equality-system} implies
$D\m\in L^\infty$. We take the Lipschitz representative of $\m$ and set
\[
U=\{x\in\R^2:\m(x)\ne\ein_3\}.
\]
Continuity makes $U$ open, and the degree condition $Q(\m)=-1$ ensures that it is nonempty. Differentiating $X\circ\M=\mathrm{id}_{\R^2}$ gives $DX(\n)\fJ(\n)=I_{\R^2}$ for $\n\ne\ein_3$. 
On each compact subset of $U$, the image of $\m$ lies in a compact subset of $\Sph\setminus\{\ein_3\}$. The local Sobolev chain rule and \eqref{eq:equality-system} give
\[
D(X\circ\m)=DX(\m)\fJ(\m)=I_{\R^2} \qquad\text{almost everywhere on } U.
\]
Choose a connected component $V$ of $U$. The weak derivative of
$X\circ\m-\mathrm{id}_{\R^2}$ vanishes on $V$, so $X(\m(x))=x-c$
almost everywhere on $V$ for some $c\in\R^2$. Applying $\M$ and using
continuity yields $\m(x)=\M(x-c)$ for every $x\in V$. At any boundary point
$x_0\in\p V$, continuity would give $\m(x_0)=\M(x_0-c)\ne\ein_3$,
contradicting $\p V\cap U=\varnothing$ because $V$ is a component of
the open set $U$. Thus $V=\R^2$.
Conversely, each translate $\M(\cdot-c)$ belongs to $\cA_{-1}$
and has energy $\cF(\M)$ by translation invariance.
\end{proof}
\section{The decreasing radial profile}\label{sec:radial}
We construct the reference map for the calibration by shooting. The argument yields a regular, strictly decreasing radial profile for $\kappa>0$ and $h>|A|$, together with the far-field asymptotics used in the next section.

For a radial profile $\theta:(0,\infty)\to \R$, set
\begin{equation}\label{eq:M}
\M(r,\varphi)=\bra*{\sin \theta(r)\cos \varphi, \, \sin \theta(r)\sin \varphi, \, \cos \theta(r)}.
\end{equation}
Direct substitution gives
\[
E(\M)=\pi \int_0^\infty \bra*{r(\theta')^2+\frac{\sin^2\theta}{r}+4\kappa r\theta'\sin^2\theta+2rW(\theta)}\dd r,
\]
where
\[
W(\theta)=h(1-\cos \theta) + \frac{A}{2}\sin^2 \theta.
\]
The corresponding Euler-Lagrange equation is~\eqref{eq:radial-ode}. For a smooth finite-energy radial map with $\theta(0)=\pi$ and $\theta(\infty)=0$, the degree is
\begin{equation}\label{eq:M_degree}
Q(\M)=\frac{1}{2}\int_0^\infty \theta' \sin \theta \dif r = \frac{1}{2}\pra*{\cos \theta(0)-\cos \theta(\infty)}=-1.
\end{equation}
\begin{proposition}\label{prop:shooting}
 Suppose $\kappa>0$ and $h>|A|$. There is a regular solution of~\eqref{eq:radial-ode} satisfying~\eqref{eq:radial-bc}. At the center it has the expansion
\begin{equation}\label{eq:center}
\theta(r)=\pi-\alpha r + \bra*{\frac{\alpha^3}{12}-\frac{\kappa\alpha^2}{4}+\frac{(h-A)\alpha}{8}}r^3+O(r^5), \quad\text{for some } \alpha>0.
\end{equation}
With $H=\sqrt{h+A}$, its tail satisfies
\begin{equation}\label{eq:tail_derivative}
-\frac{\theta'(r)}{\theta(r)}=H+\frac{1}{2r}+O(r^{-2}) \quad\text{as}\quad r\to \infty,
\end{equation}
and for some $K>0$,
\begin{equation}\label{eq:tail_K}
\theta(r)=Kr^{-\frac{1}{2}}e^{-Hr}\bra*{1+O(r^{-1})} \quad\text{as}\quad r\to \infty.
\end{equation}
The map \eqref{eq:M} is smooth on the plane and belongs to $\cA_{-1}$.
\end{proposition}
To seek a regular decreasing profile, we use the center slope $\alpha=-\theta'(0)>0$ as the shooting parameter and set $z=-r\theta'$. For $r>0$, equation~\eqref{eq:radial-ode} becomes
\begin{equation}\label{eq:shoot-system}
 \theta'=-\frac{z}{r},\qquad
 z'=-\frac{\sin\theta\cos\theta}{r}+2\kappa\sin^2\theta -r\sin\theta(h+A\cos\theta).
\end{equation}
\subsection{The singular initial-value problem}
We call a trajectory center-regular if $\theta$ extends as a $C^1$ function to $r=0$ with $\theta(0)=\pi$.
\begin{lemma}\label{lem:center-trajectories}
For each $\alpha>0$, there is a unique center-regular local solution $(\theta_\alpha,z_\alpha)$ of \eqref{eq:shoot-system} with $\theta_\alpha(0)=\pi$ and $\theta_\alpha'(0)=-\alpha$. 
For every finite $r>0$, the map $\alpha \mapsto (\theta_\alpha,z_\alpha)$ is continuous with values in $C^1([0,r];\R^2)$.
The solution satisfies the center expansion \eqref{eq:center}, and its associated radial map extends smoothly across the origin.
\end{lemma}
\begin{proof}
To treat the singularity of \eqref{eq:shoot-system} at $r=0$, write $u=\pi-\theta$ and
\begin{equation*}
 \begin{split}
 & L_0u =(ru')'-u/r,\\
 & F(r,u) =\frac{\sin u\cos u-u}{r}+2\kappa\sin^2u -hr\sin u+Ar\sin u\cos u.
 \end{split}
\end{equation*}
Then equation~\eqref{eq:radial-ode} becomes $L_0u=F(r,u)$. Setting $u=rv$ gives $L_0(rv)=r^{-1}(r^3 v')'$. For any center-regular solution with $u'(0)=\alpha$, we have $v(0)=\alpha$ and $r^3v'=r^2u'-ru\to0$ as $r\downarrow0$. Integrating twice gives the Volterra equation
\begin{equation}\label{eq:Volterra}
u =T_\alpha u := \alpha r+\cK F(\cdot,u), \qquad (\cK f)(r):=r\int_0^r s^{-3} \int_0^s \tau f(\tau) \dd \tau \dd s.
\end{equation}
We now show that $T_\alpha$ is a contraction. For $R>0$ define the space
\begin{equation*}
\begin{split}
& X_R := \set*{u \in C^1([0,R]): \, u(0)=0, \, \norm{u}_{R}<\infty},\\
& \norm{u}_{R}:=\sup_{0<r \le R} \frac{|u(r)|}{r} + \sup_{0 \le r \le R} |u'(r)|.
\end{split}
\end{equation*}
Consider the closed ball $B_{\alpha}:=\set*{u \in X_R, \,\norm{u-\alpha r}_R \le 1}$ and put $M_\alpha=\alpha+1$. Since $|u(r)|\le M_\alpha r$ on the ball $B_\alpha$, elementary estimates
\[
\abs*{\sin u \cos u - u}\le \frac{2|u|^3}{3}, \quad (\sin u \cos u-u)_u = -2\sin^2 u, \quad |(\sin^2 u)_u|\le 2|u|
\]
and the mean value theorem give, for any $u \in B_\alpha$,
\[
\abs*{F(r,u)} \le C_0 r^2, \qquad C_0 = \frac{2}{3}M_\alpha^3+2\kappa M_\alpha^2+(h+|A|)M_\alpha,
\]
and for any $u_1,u_2 \in B_\alpha$,
\[
\abs*{F(r,u_1)-F(r,u_2)}\le C_1 r \abs*{u_1-u_2}, \qquad C_1= 2M_\alpha^2+4\kappa M_\alpha+h+|A|.
\]
For $\abs*{f(r)}\le Cr^2$, integration in~\eqref{eq:Volterra} gives
\begin{equation}\label{eq:Volterra_K}
\abs*{\cK f} \le \frac{C}{8}r^3, \qquad \abs*{(\cK f)'}\le \frac{3C}{8} r^2, \qquad \norm{\cK f}_R \le \frac{C}{2}R^2.
\end{equation}
Choose $R>0$ so that $C_0R^2/2\le1$ and $q:=C_1R^2/2<1$. By~\eqref{eq:Volterra_K}, $T_\alpha$ preserves $B_\alpha$ and is a contraction with constant $q$.
Its fixed point satisfies $u'(0)=\alpha$ by~\eqref{eq:Volterra_K}. Conversely, any center-regular solution with this slope belongs to $B_\alpha$ after shrinking $R$, since $\norm{u-\alpha r}_R\to0$ as $R\downarrow0$. The contraction therefore gives uniqueness near the center. Let $(\theta_\alpha,z_\alpha)$ denote the resulting trajectory.

On a compact interval of slopes, all fixed points satisfy $|u_\alpha(r)|\le Mr$ for a common $M$. The mean value estimate for $F$ uses only this bound and therefore applies to any pair of fixed points. Choosing $R$ and $q$ uniformly using $M$ and subtracting the two Volterra equations gives
\[
\norm{u_{\alpha_1} - u_{\alpha_2}}_R \le \frac{2}{1-q} \abs*{\alpha_1 - \alpha_2}.
\]
This estimate controls $\theta_\alpha$, $\theta_\alpha'$, and $z_\alpha=ru_\alpha'=-r\theta'_\alpha$. It also controls $z_\alpha'$ through \eqref{eq:shoot-system}. Hence $(\theta_\alpha,z_\alpha)$ depends continuously on $\alpha$ in $C^1([0,R])$.
On every finite interval away from the origin, $z_\alpha'$ is bounded and $\theta_\alpha'=-z_\alpha/r$. Standard ODE continuation and continuous dependence extend the solution to every finite interval, with continuous dependence on $\alpha$ in the $C^1$ topology.

The Volterra equation~\eqref{eq:Volterra} first gives $u=\alpha r+O(r^3)$. Expanding $F(r,u)$ and substituting back yields \eqref{eq:center}.
For smoothness at the center, set $v=u/r$ and
\[
\Psi(r,\xi)=\frac{F(r,r\xi)}{r^2}.
\]
The function $\Psi$ extends smoothly and evenly across $r=0$, and the Volterra equation gives
\[
v'(r)=r\int_0^1\tau^3 \Psi(r\tau,v(r\tau))\dd\tau.
\]
Since $v'(0)=0$, the even extension is $C^1$, and the identity bootstraps it to a smooth even function. Thus $u=rv$ is smooth and odd. The functions $\sin u(r)/r$ and $-\cos u(r)$ are smooth and even, so the Cartesian map in~\eqref{eq:M} extends smoothly across the origin.
The contraction argument for $T_\alpha$ also applies at $\alpha=0$, where its unique fixed point is $u=0$, giving $\theta\equiv\pi$.
Thus every nonconstant decreasing regular profile has $\alpha>0$.
\end{proof}

Small and large center slopes $\alpha$ give rise to two different exit behaviors: a turn while $\theta$ is still positive and a crossing of zero, respectively. These two behaviors will be used to select a global decreasing profile.
\subsection{Small-slope undershooting}\label{sec:undershoot}
We call a trajectory an \emph{undershoot} if its first exit from $0<\theta<\pi$, $z>0$ occurs at a finite radius $r_*$ with $\theta(r_*)>0$, $z(r_*)=0$, and $z'(r_*)<0$. Equivalently, $\theta'(r_*)=0$ and $\theta''(r_*)>0$, so the profile turns upward while it is still positive.

For $\alpha>0$, the normalized function $w_\alpha=(\pi-\theta_\alpha)/\alpha$ satisfies $L_0w_\alpha=F_\alpha(r,w_\alpha)$ with $w_\alpha(0)=0$ and $w_\alpha'(0)=1$, where
\[
F_\alpha(r,w)=\frac{F(r,\alpha w)}{\alpha}, \qquad F_0(r,w)=-(h-A)rw.
\]
The trigonometric estimates from the center construction applied to $\cR_\alpha=F_\alpha-F_0$ give
\[
\abs*{\cR_\alpha(r,w)}\le C\alpha r^2, \qquad
\abs*{\partial_w\cR_\alpha(r,w)}\le C\alpha r
\]
uniformly for $0<r\le1$, $|w|\le2r$, and $0<\alpha\le1$. By~\eqref{eq:Volterra_K}, for small $\rho>0$ the map $S_\alpha w=r+\cK F_\alpha(\cdot,w)$ preserves the ball $\norm{w-r}_\rho\le1$ in $X_\rho$ and contracts it with a common constant $q<1$ for $\alpha\in[0,1]$. For $\alpha>0$ its fixed point is $w_\alpha$ by center uniqueness. Comparing with the fixed point $w_0$ of $S_0$ gives
\[
\norm{w_\alpha-w_0}_\rho \le q\norm{w_\alpha-w_0}_\rho+C\alpha\rho^2.
\]
On each fixed interval $[\rho,R]$, the normalized equation is regular. There $F_\alpha$ differs from $F_0$ by $O_R(\alpha)$ and is uniformly Lipschitz on bounded sets of $w$. ODE stability and the center estimate thus give
\begin{equation}\label{eq:w_convergence}
\norm{w_\alpha-w_0}_{C^1([0,R])}\le C_R\alpha
\end{equation}
for sufficiently small $\alpha$.
The limiting equation can be written as
\[
w_0''+\frac{1}{r} w_0' + \bra*{h-A-\frac{1}{r^2}} w_0=0,
\]
which is the order-one Bessel equation after rescaling $r$ by $\sqrt{h-A}$. The center condition $w_0(0)=0$ excludes the singular branch $Y_1$, while $w_0'(0)=1$ fixes the normalization. Thus
\[
w_0(r)=\frac{2}{\sqrt{h-A}}J_1(\sqrt{h-A} \, r).
\]
Let $r_1$ be the first positive zero of $w_0$ and $r_m\in(0,r_1)$ its first critical point, which exists because $w_0'(0)=1$ and $w_0(r_1)=0$. On $(0,r_1)$ we have $w_0>0$. Since $(rw_0')'=(r^{-1}-(h-A)r)w_0$, the function $rw_0'$ is positive up to $(h-A)^{-1/2}$ as long as $w_0>0$. Hence $r_m>(h-A)^{-1/2}$ and $w_0''(r_m)<0$. Choose $0<r_-<r_m<r_+<r_1$ with $w_0'>0$ on $[0,r_-]$ and $w_0'(r_+)<0$.

By~\eqref{eq:w_convergence}, these derivative signs persist for small $\alpha$. Thus $w_\alpha>0$ on $(0,r_-]$ since $w_\alpha(0)=0$, and $w_\alpha>0$ on $[r_-,r_+]$ by uniform convergence to $w_0>0$ there. Boundedness of $w_\alpha$ gives $0<\theta_\alpha=\pi-\alpha w_\alpha<\pi$ on $(0,r_+]$ for small $\alpha$. The first zero of $w_\alpha'$ occurs in $(r_-,r_+)$ and gives the first contact with $z_\alpha=0$, since $z_\alpha=\alpha r w_\alpha'$.

The following sign lemma will also be used in Section~\ref{sec:barriers}.
\begin{lemma}\label{lem:lower-crossing-zeros}
Assume $\kappa>0$ and $h>|A|$. For $\nu\in[0,1]$, define
\begin{align*}
& L(\nu)=h+A(1-\nu^2),\\
& P(\nu)=h(3\nu^2-1)+A(1-\nu^2)(4\nu^2-1),\\
& g_\sigma(\nu)=P(\nu)-4\sigma\kappa\nu(1-\nu^2)\sqrt{L(\nu)}, \quad \sigma \in \set*{-1,1}.
\end{align*}
Each $g_\sigma$ has exactly one zero $\nu_\sigma \in (0,1)$, is negative on $[0,\nu_\sigma)$ and positive on $(\nu_\sigma,1]$, and satisfies
\[
(g_\sigma)_\nu(\nu_\sigma)>0.
\]
\end{lemma}
\begin{proof}
Set
\[
 f(\nu)=-\frac{P(\nu)}{\nu(1-\nu^2)\sqrt{L(\nu)}}.
\]
Direct differentiation gives
\[
f_\nu(\nu)= -\frac{(h+A)^2}{\nu^2 L^{3/2}}-\frac{2(L+h\nu^2)}{\sqrt{L}(1-\nu^2)^2}<0 \quad\text{for all}\quad \nu \in (0,1),
\]
since $h>0$ and $L>0$. The endpoint values $P(0)=-(h+A)$ and $P(1)=2h$ give $f(0_+)=+\infty$ and $f(1_-)=-\infty$. Hence $f(\nu)=-4\sigma\kappa$ has exactly one solution $\nu_\sigma \in (0,1)$, which is the unique zero of $g_\sigma$. Moreover, 
\[
g_\sigma(\nu)=-\nu(1-\nu^2)\sqrt{L(\nu)}\bra*{f(\nu)+4\sigma\kappa},
\]
gives the asserted signs and
\[
(g_\sigma)_\nu(\nu_\sigma)=-\nu_\sigma(1-\nu_\sigma^2)\sqrt{L(\nu_\sigma)} f_\nu(\nu_\sigma) >0.
\]
\end{proof}

\begin{lemma}\label{lem:transverse}
Assume $\kappa>0$ and $h>|A|$. Let $r_*>0$ be the first positive zero of $z$ along a center-regular trajectory with slope $\alpha>0$. If $\theta(r_*)>0$, then $z'(r_*)<0$.
\end{lemma}

\begin{proof}
By the choice of $r_*$, we have $z(r)>0$, $\theta'(r)<0$ and $\theta(r) \in (0,\pi)$ for any $r\in (0,r_*)$, and hence $z'(r_*)\le 0$. Suppose for contradiction that the zero-slope contact $z(r_*)=0$ is a tangency, that is $z'(r_*)= 0$. Taylor expansion at $r_*$, together with $z(r)>0$ for $r<r_*$, gives $z''(r_*)\ge 0$.
Set
\[
\nu = \sin \frac{\theta}{2}, \qquad Q(\nu)=h+A(1-2\nu^2)>0,
\]
and use subscript $_*$ for evaluation at $r=r_*$. Differentiating the $z$-equation in~\eqref{eq:shoot-system} and using $\theta'(r_*)=0$ yield
\[
z''(r_*)=\sin \theta_*\bra*{\frac{\cos \theta_*}{r_*^2}-Q_*}.
\]
Thus $z'(r_*)=0$ and $z''(r_*)\ge 0$ imply
\begin{equation}\label{eq:grazing}
Q_* r^2_*-2\kappa r_* \sin \theta_* + \cos \theta_*=0, \qquad Q_* r^2_* \le \cos \theta_*.
\end{equation}
In particular, since $Q>0$, we have $\cos \theta_*>0$ and $0<\nu_*<1/\sqrt{2}$.

We derive opposite signs for $g_-(\nu_*)$ from the trajectory before contact and from the tangency equations. Set
\[
V=-\theta'-\frac{\sin\theta}{r}, \qquad S=\sqrt{2W(\theta)}.
\]
The tangency conditions imply $V_*<-S_*$, as shown below. We use $F=S^2-V^2$ to locate the last contact with $V=-S$ before $r_*$ and examine the sign of $F'$ there.
The center expansion~\eqref{eq:center} gives $V(r) \to 0$ and $W(\theta(r)) \to 2h$, so $F(r) \to 4h>0$ as $r \downarrow 0$. The identity
\[
\sin^2\theta \, Q(\nu) - 2\cos \theta\, W(\theta)=h(1-\cos \theta)^2>0
\]
and~\eqref{eq:grazing} give
\[
V^2_*=\frac{\sin^2 \theta_*}{r^2_*} \ge \frac{\sin^2 \theta_* Q_*}{\cos \theta_*} > 2W(\theta_*)=S^2_*.
\]
Since $V_*=-\sin \theta_*/r_*<0$, we have $V_*<-S_*$ and $F(r_*)<0$.
With $L,P,g_-$ as in Lemma~\ref{lem:lower-crossing-zeros}, we have
\[
P=2\nu^2L-(1-\nu^2)Q, \quad W=2\nu^2 L, \quad S=2\nu\sqrt{L} \quad\text{and}\quad W_\theta=\sin \theta\, Q(\nu).
\]
The radial equation~\eqref{eq:radial-ode} gives
\[
r V' =-(1-\cos \theta)V+2\kappa \sin^2\theta-rW_\theta.
\]
Consequently, 
\[
rF'=-2\sin \theta\,W_\theta+2(1-\cos\theta)V^2-4\kappa V \sin^2\theta.
\]
Substituting $V=-S$ and $W_\theta=\sin\theta\,Q(\nu)$ gives
\begin{equation}\label{eq:no-grazing-cross}
rF'=8\nu^2g_-(\nu)\qquad\text{on }V=-S.
\end{equation}
Thus, on the boundary $V=-S$, the sign of $F'$ is the sign of $g_-(\nu)$. Let $r_1 \in (0,r_*)$ so that $F(r_1)=0$ and $F(r)<0$ for any $r\in (r_1,r_*]$. Hence $F'(r_1)\le 0$. Since $V_*<-S_*$, the continuity forces $V<-S$ throughout the interval $(r_1,r_*]$, and hence $V(r_1)=-S(r_1)$. Writing $\nu_1=\nu(r_1)$,~\eqref{eq:no-grazing-cross} and $\nu'=\frac{1}{2}\cos(\theta/2)\theta'<0$ give
\[
g_-(\nu_1)\le 0 \quad\text{and}\quad \nu_*<\nu_1.
\]
By Lemma~\ref{lem:lower-crossing-zeros}, $g_-$ has a unique zero and $g_-(\nu_*)<0$.

At the tangency point $r_*$, the equality in~\eqref{eq:grazing} gives
\begin{equation}\label{eq:grazing-g-2}
g_-(\nu_*)-P_*=\sqrt{(1-\nu^2_*)L_*}\bra*{Q_*r_*+\frac{1-2\nu^2_*}{r_*}}\ge 2\sqrt{(1-\nu^2_*)(1-2\nu^2_*)L_*Q_*}>0.
\end{equation}
If $P_*\ge 0$, then $g_-(\nu_*)>0$ immediately. If $P_* <0$, its definition gives $2\nu^2_*L_*<(1-\nu^2_*)Q_*$, and hence
\begin{align*}
& 4(1-\nu^2_*)(1-2\nu^2_*)L_*Q_*-P^2_* =(1-\nu^2_*)^2 Q_*(4L_*-Q_*)-4\nu^4_*L^2_*\\
>&(1-\nu^2_*)^2 Q_*(4L_*-Q_*)-(1-\nu^2_*)^2Q^2_*=2(1-\nu^2_*)^2 Q_*(2L_*-Q_*)\\
=&2(1-\nu^2_*)^2 Q_*(h+A)>0.
\end{align*}
The positive lower bound in~\eqref{eq:grazing-g-2} exceeds $-P_*$, so $g_-(\nu_*)>0$ also in this case. This contradicts the inequality $g_-(\nu_*)<0$ established above. Hence the assumed tangency is impossible, and $z'(r_*)<0$.
\end{proof}

Consequently, every sufficiently small positive center slope produces an undershoot.
\subsection{Large-slope overshooting}
We call a trajectory an \emph{overshoot} if its first exit from $0<\theta<\pi$, $z>0$ occurs through $\theta=0$ with $z>0$ at a finite radius. Thus the profile $\theta$ crosses zero while still decreasing.

For a trajectory with center slope $\alpha>0$, introduce the variable $s=\alpha r$ and write $U(s)=\pi-\theta(s/\alpha)$. The scaled equation is
\[
(sU_s)_s=\frac{\sin U \cos U}{s}+\frac{2\kappa}{\alpha}\sin^2 U - \frac{s}{\alpha^2} \sin U(h-A\cos U), \qquad U_s(0)=1.
\]
For large $\alpha$, the fixed potential enters as a small perturbation. We therefore begin with the zero-potential system $h=A=0$, obtain quantitative bounds on its zero crossing, and then control the effect of the potential to show the full trajectory also overshoots.

\paragraph{Quantitative bounds for the zero-potential crossing}

Keep $\kappa$ fixed and let $U_\alpha^0(s)=\pi-\theta_\alpha^0(s/\alpha)$ be the zero-potential trajectory. The center construction and continuation in Lemma~\ref{lem:center-trajectories} apply also to $h=A=0$. The scaled equation gives
\begin{equation}\label{eq:D-identity}
D=(s(U_\alpha^0)_s)^2-\sin^2U_\alpha^0,\qquad
D_s=\frac{4\kappa}{\alpha}s(U_\alpha^0)_s\sin^2U_\alpha^0.
\end{equation}
As $\alpha\to\infty$, the scaled zero-potential equation tends to $(sU_s)_s=\sin U\cos U/s$. Its center-regular solution with $U(0)=0$ and $U_s(0)=1$ is the Belavin--Polyakov profile $U_{\mathrm{BP}}(s)=2\arctan(s/2)$. On the center ball, where $|U(s)|\le 2s$, the perturbation $(2\kappa/\alpha)\sin^2 U$ is $O((\kappa/\alpha)s^2)$. The Volterra and compact-interval stability argument used in Section~\ref{sec:undershoot} gives, for each fixed $S$ and all sufficiently large $\alpha$,
\begin{equation}\label{eq:BP-limit}
\norm{U_\alpha^0-U_{\mathrm{BP}}}_{C^1([0,S])}\le C_S\frac{\kappa}{\alpha}.
\end{equation}
Since $(U_{\mathrm{BP}})_s\ge1/2$ on $[0,2]$ and $U_{\mathrm{BP}}(2)=\pi/2$, we have $(U_\alpha^0)_s\ge1/4$ there and $U_\alpha^0(2)<\pi$ for large $\alpha$. Integrating~\eqref{eq:D-identity} and using the convergence on $[1,2]$ gives
\begin{equation}\label{eq:D-seed}
D(2)=\frac{4\kappa}{\alpha}\int_0^2 s(U_\alpha^0)_s\sin^2U_\alpha^0\,\dd s
\ge c_0\frac{\kappa}{\alpha}
\end{equation}
for some $c_0>0$ independent of $\alpha$. For $s\ge2$, the identity~\eqref{eq:D-identity} gives $D(s)\ge D(2)$ as long as $(U_\alpha^0)_s>0$. This excludes a first zero of $(U_\alpha^0)_s$ and yields $s(U_\alpha^0)_s\ge\sqrt{c_0\kappa/\alpha}$ also beyond the crossing of $U_\alpha^0=\pi$.

Let $s_0>2$ be the first radius where $U_\alpha^0=\pi$, with $s_0=\infty$ if there is no such radius. Before $s_0$, we have $D>0$ by~\eqref{eq:D-identity} and $D(0)=0$, so $s(U_\alpha^0)_s>\sin U_\alpha^0$. Thus
\[
\frac{\dif}{\dif s}\log\bra*{s\cot\frac{U_\alpha^0}{2}}
=\frac{1-s(U_\alpha^0)_s/\sin U_\alpha^0}{s}<0.
\]
Since $s\cot(U_\alpha^0/2)\to2$ as $s\downarrow0$, this gives
\begin{equation}\label{eq:U_alpha}
\pi-U_\alpha^0(s)\le 2\arctan(2/s)\le 4/s\qquad(0<s<s_0).
\end{equation}
Set $s_1=\sqrt{\alpha/\kappa}$ and choose $L_{\rm ex}>1$ with $\sqrt{c_0}\log L_{\rm ex}>4$. If $s_0\le s_1$, then $s_0<L_{\rm ex}s_1$. Otherwise $\pi-U_\alpha^0(s_1)\le4\sqrt{\kappa/\alpha}$. If $s_0\ge L_{\rm ex}s_1$, integrating the speed bound from $s_1$ to $L_{\rm ex}s_1$ gives
\[
\pi-U_\alpha^0(L_{\rm ex}s_1)
\le\sqrt{\frac{\kappa}{\alpha}}\bra*{4-\sqrt{c_0}\log L_{\rm ex}}<0,
\]
contrary to $U_\alpha^0\le\pi$ up to $s_0$. Thus $s_0<L_{\rm ex}\sqrt{\alpha/\kappa}$.

\paragraph{Error bounds on a growing interval}

We compare the two trajectories on $0<s\le2L_{\rm ex}\sqrt{\alpha/\kappa}$. In the variable $t=\log s$, write $\Theta(t)=\pi-U(e^t)$ and $Z(t)=-\Theta_t=sU_s$. Use $(\Theta_\alpha,Z_\alpha)$ for the full potential and $(\Theta_\alpha^0,Z_\alpha^0)$ for zero potential. For the full trajectory,
\[
(\Theta_\alpha)_{tt}
=\sin\Theta_\alpha\cos\Theta_\alpha-\frac{2\kappa s}{\alpha}\sin^2\Theta_\alpha
+\frac{s^2}{\alpha^2}\sin\Theta_\alpha(h+A\cos\Theta_\alpha).
\]
The zero-potential equation omits the last term. For $\delta=\Theta_\alpha-\Theta_\alpha^0$, the mean value formula on the first two terms gives $\delta_{tt}=C(t)\delta+F(t)$, where $F(t)$ is the full potential term and
\begin{equation}\label{eq:difference_C_F}
|C(t)|\le 1+\frac{2\kappa s}{\alpha},\qquad
|F(t)|\le \frac{(h+|A|)s^2}{\alpha^2}.
\end{equation}
The first-order system for $Y=(\delta,\delta_t)^{\fT}$ has coefficient matrix with rows $(0,1)$ and $(C(t),0)$. On the comparison interval, $2\kappa s/\alpha\le4L_{\rm ex}\sqrt{\kappa/\alpha}$, so the matrix has Euclidean operator norm at most $M=3/2$ for sufficiently large $\alpha$. The homogeneous propagator therefore has norm bounded by $e^{M(t-\tau)}$ for $\tau\le t$. The center expansions give $|Y(T)|=O_\alpha(e^{3T})$ as $T\to-\infty$, so the propagated initial term tends to zero. Since $M<2$, Duhamel's formula gives
\begin{equation}\label{eq:Duhamel}
|\delta(t)|+|\delta_t(t)|
\le\frac{C(h+|A|)}{\alpha^2}\int_{-\infty}^t e^{M(t-\tau)}e^{2\tau}\,\dd\tau
\le c_1\frac{(h+|A|)s^2}{\alpha^2}.
\end{equation}
The last inequality follows from the integral being $e^{2t}/(2-M)$. In particular, the error on the comparison interval is $O((h+|A|)/(\kappa\alpha))=o(\sqrt{\kappa/\alpha})$ as $\alpha\to\infty$ with the parameters fixed.

\paragraph{Overshooting for the full-potential trajectory}

Let $t_0=\log s_0$ and $t_1=t_0+\log 2$. The corresponding radius $s=2s_0$ lies within the comparison interval. The lower speed bound from~\eqref{eq:D-seed} gives
\[
\Theta_\alpha^0(t_1)=-\int_{t_0}^{t_1}Z_\alpha^0(t)\,\dd t
\le-(\log 2)\sqrt{c_0\kappa/\alpha}.
\]
Thus $\Theta_\alpha(t_1)<0$ for large $\alpha$ by~\eqref{eq:Duhamel}. To exclude an earlier exit through zero speed, use $Z_\alpha^0\ge s/4$ for $0<s\le2$ and the local error bound $O(s^2/\alpha^2)$ to obtain $Z_\alpha\ge s/8>0$ there. For $2\le s\le2s_0$, use $Z_\alpha^0\ge\sqrt{c_0\kappa/\alpha}$ and the uniform error $o(\sqrt{\kappa/\alpha})$ to obtain $Z_\alpha>0$. The full profile therefore first exits through $\theta=0$ while still decreasing, at a physical radius smaller than $2L_{\rm ex}/\sqrt{\kappa\alpha}$.

\subsection{The connecting trajectory and its tail}

\begin{proof}[Proof of Proposition~\ref{prop:shooting}]
For each center-regular trajectory, define its first positive exit time from the region $0<\theta<\pi$, $z>0$ by
\[
T(\alpha)=\sup \set*{R>0: \, 0<\theta_\alpha<\pi, \, \theta_\alpha' <0 \text{ on } (0,R)}.
\]
By Lemma~\ref{lem:center-trajectories}, the trajectory exists at every finite radius. If $T(\alpha)<\infty$, the strict decrease of $\theta_\alpha$ rules out an exit through $\theta=\pi$. A finite exit cannot occur at $(\theta,z)=(0,0)$, since ODE uniqueness at a positive radius would force the trajectory to be identically zero, contradicting $\theta(0)=\pi$.  A contact with $z=0$ at $0<\theta<\pi$ is transverse by Lemma~\ref{lem:transverse}, with $z'<0$.
An exit through $\theta=0$, $z>0$, is transverse because $\theta'=-z/r<0$. Thus every finite exit is either an undershoot or an overshoot.

Fix a slope $\alpha_0$ with finite exit radius $T_0$. For any $t\in(0,T_0)$, both $\theta_{\alpha_0}$ and $-\theta'_{\alpha_0}$ have positive minima on $[0,t]$. By $C^1$ continuous dependence and $\theta_\alpha(0)=\pi$, nearby trajectories satisfy the strict region inequalities on $(0,t]$. These inequalities persist just beyond $t$, so $T(\alpha)>t$.

The exiting coordinate is $\theta$ for an overshoot and $z$ for an undershoot. By transversality it is negative at $T_0+\delta$ for small $\delta>0$, while the other region inequalities hold strictly on $[T_0-\delta,T_0+\delta]$. By continuous dependence, the same inequalities and endpoint sign hold for nearby slopes. Taking $t=T_0-\delta$ in the preceding argument forces their first exit to occur through the same boundary within this interval. Thus both sets of slopes are open.

The two sets are disjoint and nonempty by the small- and large-slope results. They cannot cover the connected interval $(0,\infty)$, so there exists a slope with $T(\alpha)=\infty$. Its profile decreases to some $\theta_*\in[0,\pi)$. If $\theta_*>0$, then $\sin\theta$ is bounded below by a positive constant for large $r$. The equation for $z$ in~\eqref{eq:shoot-system} and $h>|A|$ give $z'\le-cr$ for large $r$ and some $c>0$. This would force $z$ to become negative. Hence $\theta_*=0$.

For the tail estimate, set $H=\sqrt{h+A}>0$ and write the equation~\eqref{eq:radial-ode} as
\begin{equation}\label{eq:radial}
\theta''+\frac{1}{r}\theta'-\bra*{H^2+\frac{1}{r^2}}\theta=\cN(r,\theta),
\end{equation}
where
\[
\cN(r,\theta)=\frac{\sin\theta\cos\theta-\theta}{r^2}-\frac{2\kappa\sin^2\theta}{r}+h(\sin\theta-\theta)+A(\sin\theta\cos\theta-\theta).
\]
Along the profile, $\cN(r,\theta)=O(\theta^2/r+\theta^3)$ as $r\to\infty$. Put $\eta=\cN(r,\theta)/\theta=O(\theta/r+\theta^2)$, so that $\eta \to 0$ as $r\to\infty$. Fix $H_0\in(0,H)$, and choose $R$ so large that $Q(r):=H^2+r^{-2}+\eta(r)\ge H_0^2$ for $r \ge R$. The operator
\[
\cL=\p_{rr}+\frac{1}{r}\p_r-Q(r)
\]
satisfies $\cL\theta=0$. Define $\phi(r)=e^{-H_0(r-R)}$. Then
\[
\cL\phi=\bra*{H_0^2-\frac{H_0}{r}-Q(r)}\phi\le 0
\qquad\text{for }r\ge R.
\]
Set $\omega(r)=\theta(r)-\theta(R)\phi(r)$. Then $\omega(R)=0$, $\omega(r)\to0$ as $r\to\infty$, and $\cL\omega\ge0$ on $[R,\infty)$. If $\omega$ were positive somewhere, it would attain a positive maximum at an interior point $r_0>R$.
At this point,
\[
\cL\omega(r_0) =\omega''(r_0)-Q(r_0)\omega(r_0)<0,
\]
because $\omega''(r_0)\le0$ and $Q(r_0)\ge H_0^2>0$. This is a contradiction. Hence $\omega\le0$, and therefore
\begin{equation}\label{eq:preliminary-decay}
0<\theta(r)\le\theta(R)e^{-H_0(r-R)} \qquad\text{for }r\ge R.
\end{equation}
In particular $\eta(r)=O(e^{-H_0 r})$.
To remove the first-order term, set $y(r)=\sqrt{r}\theta(r)$. Direct substitution into~\eqref{eq:radial} gives
\[
y''=(H^2+g(r))y, \qquad g(r)=\frac{3}{4r^2}+\eta(r)=O(r^{-2}).
\]
Here $y>0$ and $y \to 0$ by~\eqref{eq:preliminary-decay}. Increasing $R$ if necessary, we have $H^2/4 \le H^2+g(r)\le 4H^2$ for $r \ge R$. On this interval $y''>0$, so $y\to0$ implies $y'<0$ and $y'\uparrow0$. Moreover,
\[
(y'^2-4H^2y^2)'=2y'y(g-3H^2) \ge 0.
\]
Since $y,y'\to0$, this yields $y'^2\le4H^2y^2$, so $q=-y'/y$ satisfies $0<q\le2H$. Integrate $(q-H)'=(q+H)(q-H)-g$ over $[r,T]$ with an integrating factor. The terminal term has modulus at most $He^{-H(T-r)}$, so letting $T\to\infty$ gives
\[
q(r)-H=\int_r^\infty e^{-\int_r^s(q(t)+H)\,\dif t} g(s)\,\dif s.
\]
Since $g(s)=O(s^{-2})$ and $q+H\ge H$, we obtain
\[
|q(r)-H|\le C\int_r^\infty e^{-H(s-r)}s^{-2}\,\dif s\le \frac{C}{Hr^2}.
\]
Now $-\theta'/\theta=q+1/(2r)$, which proves~\eqref{eq:tail_derivative}.

Let $F(r)=e^{Hr}\sqrt r\,\theta(r)>0$. Since $(\log F)'=H-q=O(r^{-2})$, $\log F$ converges to a finite limit $\log K$ with $K>0$, and $\log(F(r)/K)=O(r^{-1})$. This proves~\eqref{eq:tail_K}.

Finally, we consider the regularity of the corresponding $\M$ given in~\eqref{eq:M}. The logarithmic-derivative estimate~\eqref{eq:tail_derivative} also gives $|\theta'| = O(\theta)$ on the tail. The radial equation~\eqref{eq:radial} then yields
\[
\theta''=\bra*{H^2+\frac{1}{r^2}+\eta}\theta-\frac{1}{r}\theta'=O(\theta).
\]
The polar-coordinate derivative formulas give
\[
|D\M| \le C\bra*{|\theta'|+\frac{|\sin \theta|}{r}} = O(\theta) \quad\text{and}\quad
|D^2\M| \le C\bra*{|\theta''| + |\theta'|^2+\frac{|\theta'|}{r}+\frac{|\sin \theta|}{r^2}} =O(\theta).
\]
Moreover, $|\M-\ein_3|=O(\theta)$ as $r\to\infty$. Combining these tail estimates with the smooth extension of $\M$ through the origin established in Lemma~\ref{lem:center-trajectories}, we obtain
\begin{equation}\label{eq:M_regularity}
\M-\ein_3\in H^2(\R^2;\R^3),
\qquad
D\M,D^2\M\in L^\infty(\R^2).
\end{equation}
In particular, $\M\in\cA$, and \eqref{eq:M_degree} gives $Q(\M)=-1$. Thus $\M\in\cA_{-1}$, completing the proof.
\end{proof}
\section{The two global singular-value barriers}\label{sec:barriers}
Assume $\kappa>0$ and $h>|A|$. Throughout this section, let $\theta$ be any regular solution of~\eqref{eq:radial-ode} satisfying~\eqref{eq:radial-bc}. Write $\alpha=-\theta'(0)>0$ and $H=\sqrt{h+A}$. The center and tail arguments of Section~\ref{sec:radial} apply to every such solution. The radial and angular singular values of $D\M$ for $r>0$ are
\[
a(r)=-\theta'(r)>0,\qquad
b(r)=\frac{\sin\theta(r)}r>0.
\]
For the calibration coefficient $\beta$ constructed in Section~\ref{sec:construction}, the lower barrier $(a-b)^2<2W(\theta)$ ensures $1+\beta>0$, and the upper barrier $2W(\theta)<(a+b)^2$ ensures $1-\beta>0$, by~\eqref{eq:beta-gap}. We prove the lower barrier for every such profile, characterize the upper barrier in terms of the center slope, and give a parameter condition that ensures the strict slope bound needed at the south pole.

With $W(\theta)=h(1-\cos\theta)+\frac A2\sin^2\theta$ as before, we set
\[
U=a+b,\qquad V=a-b,\qquad S=\sqrt{2W(\theta)}.
\]
The definitions of $V$ and $S$ agree with those in the proof of Lemma~\ref{lem:transverse}. The radial equation gives
\begin{equation}\label{eq:U}
rU'=(1+\cos\theta)\pra*{2\kappa(1-\cos \theta)-U} -rW_\theta(\theta),
\end{equation}
\begin{equation}\label{eq:V}
rV'=(1-\cos\theta)\pra*{2\kappa(1+\cos \theta)-V} -rW_\theta(\theta).
\end{equation}
The center expansion~\eqref{eq:center} yields
\begin{equation}\label{eq:center-signs}
 \begin{split}
& a=\alpha+O(r^2), \quad b=\alpha+O(r^2), \quad W(\theta)=2h+O(r^2),\\
& W_\theta(\theta)=\alpha(h-A)r+O(r^3), \quad
S=2\sqrt{h}+O(r^2) \quad\text{as}\quad r\to 0.
 \end{split}
\end{equation}
The tail estimates~\eqref{eq:tail_derivative} and \eqref{eq:tail_K} give
\[
a=H\theta+\frac{\theta}{2r}+O(\theta/r^2), \quad b=\frac{\theta}{r}+O(\theta^3/r), \quad S=H\theta+O(\theta^3) \quad\text{as}\quad r\to \infty.
\]
Since $\theta$ decays exponentially, it follows that
\begin{equation}\label{eq:tail-signs}
 \begin{split}
 U-S&=\frac{3\theta}{2r}+O(\theta/r^2)>0,\\
 S-V&=\frac{\theta}{2r}+O(\theta/r^2)>0,\\
 V+S&=2H\theta+O(\theta/r)>0
 \end{split}
\end{equation}
for all sufficiently large $r$.

\begin{proposition}\label{prop:lower}
 Assume $h>|A|$. Let $\theta$ be a regular solution of~\eqref{eq:radial-ode} satisfying~\eqref{eq:radial-bc}. Then
 \[
(a(r)-b(r))^2 < 2W(\theta(r)) \quad\text{for every}\quad r\in(0,\infty).
\]
\end{proposition}
\begin{proof}
Set $F=S^2-V^2$. We first exclude negative intervals of $F$ by comparing the crossing directions at their entry and exit points. Since $SS_\theta=W_\theta$ and $\theta'+V=-\sin\theta/r$, equation~\eqref{eq:V} yields
\begin{equation}\label{eq:F_whole}
rF'=\mathcal{H}(\theta,V)=2V^2(1-\cos\theta)-2\sin\theta W_\theta-4\kappa V\sin^2 \theta.
\end{equation}
Put $\nu=\sin(\theta/2)$ and use $L,g_\sigma$ and $\nu_\sigma$ from Lemma~\ref{lem:lower-crossing-zeros}. On either branch $V=\sigma S$, with $\sigma\in\{-1,1\}$, substituting $W=2\nu^2L$ and $S=2\nu\sqrt L$ into~\eqref{eq:F_whole} gives
\begin{equation}\label{eq:F_contact}
rF'=8\nu^2 g_\sigma(\nu) \quad\text{on}\quad V=\sigma S.
\end{equation}
The center expansion~\eqref{eq:center-signs} gives $F(0_+)=4h>0$, while~\eqref{eq:tail-signs} gives $F>0$ for all sufficiently large $r$.
Suppose that $\{r>0:F(r)<0\}$ is nonempty, and let $(r_1,r_2)$ be one of its connected components. Then $0<r_1<r_2<\infty$ and $F(r_1)=F(r_2)=0$.
Since $|V|>S>0$ on $(r_1,r_2)$, the sign of $V$ is fixed there, and both endpoints lie on the same branch $V=\sigma S$.
At the entry point, $F'(r_1)\le0$, so \eqref{eq:F_contact} gives $g_\sigma(\nu(r_1))\le0$ and hence $\nu(r_1)\le\nu_\sigma$.
Since $\nu$ is strictly decreasing, $\nu(r_2)<\nu(r_1)\le\nu_\sigma$. Lemma~\ref{lem:lower-crossing-zeros} therefore gives $g_\sigma(\nu(r_2))<0$.
By~\eqref{eq:F_contact}, this implies $F'(r_2)<0$, contradicting $F'(r_2)\ge0$ at the exit point. Thus $F\ge0$ on $(0,\infty)$.

It remains to exclude zeros. If $F(r_0)=0$, then $V(r_0)=\sigma S(r_0)$ for some $\sigma$. Since $F\ge0$, we have $F'(r_0)=0$ and $F''(r_0)\ge0$. By~\eqref{eq:F_contact}, the first equality gives $\nu(r_0)=\nu_\sigma$. Differentiating $F=S^2-V^2$ at $r_0$ also gives $V'(r_0)=\sigma S_\theta(\theta(r_0))\theta'(r_0)$. Thus the trajectory is tangent to the branch $V=\sigma S(\theta)$ there, and the derivative of $\mathcal H$ along the trajectory equals its derivative along that branch at $r_0$. Since $(rF')'(r_0)=r_0F''(r_0)$, differentiating~\eqref{eq:F_whole} with the branch formula~\eqref{eq:F_contact} gives
\[
r_0 F''(r_0)=8\nu_\sigma^2(g_\sigma)_\nu(\nu_\sigma)\nu'(r_0)<0.
\]
Here $(g_\sigma)_\nu(\nu_\sigma)>0$ by Lemma~\ref{lem:lower-crossing-zeros} and $\nu'(r_0)<0$. This contradicts $F''(r_0)\ge0$, so $F>0$ on $(0,\infty)$.
\end{proof}
We next characterize the upper barrier in terms of the center slope.

\begin{proposition}\label{prop:upper}
 Assume $h>|A|$. Let $\theta$ be a regular solution of~\eqref{eq:radial-ode} satisfying~\eqref{eq:radial-bc}, with center slope $\alpha=-\theta'(0)$. Then
 \[
(a(r)+b(r))^2 > 2W(\theta(r)) \quad\text{for every}\quad r>0 \qquad\Leftrightarrow\qquad \alpha \ge \sqrt{h}.
\]
\end{proposition}
\begin{proof}
Set $G=U-S$. Since $U+S>0$, the upper barrier is equivalent to $G>0$ on $(0,\infty)$. By~\eqref{eq:center-signs}, $G(0_+)=2(\alpha-\sqrt h)$, so $\alpha\ge\sqrt h$ is necessary.

For sufficiency, use~\eqref{eq:U}, $S_\theta=W_\theta/S$, and $U+\theta'=\sin\theta/r$ to obtain
\begin{equation}\label{eq:G}
G'=\bra*{\frac{W_\theta}{S}-\frac{1+\cos\theta}{r}}G+\frac{8\nu^2(1-\nu^2)}{r}\bra*{\kappa-\cR(\nu)},
\end{equation}
where $\nu=\sin(\theta/2)$, $L(\nu)=h+A(1-\nu^2)>0$, and
\[
\cR(\nu)=\frac{3h+3A-4A\nu^2}{4\nu\sqrt{L(\nu)}}, \qquad \nu \in (0,1].
\]
Since $W_\theta(\theta)/S=O(r)$ and $(1+\cos \theta)/r=O(r)$, the integrating factor
\[
I(r)=\exp\pra*{\int_0^r \bra*{\frac{1+\cos\theta(s)}{s}}-\frac{W_\theta(\theta(s))}{S(s)} \dd s}
\]
is well defined with $I(0)=1$. Multiplying~\eqref{eq:G} by $I$ gives
\begin{equation}\label{eq:IG}
(IG)'=\frac{8I\nu^2(1-\nu^2)}{r}\pra*{\kappa-\cR(\nu)}.
\end{equation}
Assume now that $\alpha\ge\sqrt h$. We have $(IG)(0_+)=G(0_+)\ge0$ and, by~\eqref{eq:tail-signs}, $I(R)G(R)>0$ for every sufficiently large $R$. Direct differentiation gives
\begin{equation}\label{eq:R}
\cR_\nu(\nu)
=-\frac{(h+A)\pra*{3(h+A)-2A\nu^2}}{4\nu^2L(\nu)^{3/2}}<0,
\qquad \cR(1)=\frac{3h-A}{4\sqrt h},
\end{equation}
where we used $h+A>0$ and $3(h+A)-2A\nu^2\ge3(h-|A|)>0$.
Since $\nu'(r)<0$ for $r>0$, the function $r\mapsto\kappa-\cR(\nu(r))$ is strictly decreasing. By~\eqref{eq:IG}, the derivative of $IG$ can change sign only from positive to negative. Fix $r_0>0$. If $\kappa-\cR(\nu(r_0))\ge0$, then $(IG)'>0$ on $(0,r_0)$, so
\[
I(r_0)G(r_0)>\lim_{r\downarrow0}I(r)G(r)\ge0.
\]
If $\kappa-\cR(\nu(r_0))<0$, then $(IG)'<0$ on $(r_0,\infty)$. Choosing a sufficiently large $R>r_0$ gives
\[
I(r_0)G(r_0)>I(R)G(R)>0.
\]
Thus $G(r_0)>0$ in either case, proving the upper barrier.
\end{proof}
Under $h>|A|$, Propositions~\ref{prop:lower} and~\ref{prop:upper} give
\[
(a-b)^2<2W(\theta)<(a+b)^2, \qquad\text{for every }r>0
\]
for every regular strictly decreasing profile $\theta$ with center slope $\alpha\ge \sqrt{h}$.

We next give an explicit parameter condition ensuring the strict center-slope bound $\alpha > \sqrt{h}$.
\begin{corollary}\label{coro:range}
Let $\theta$ be a regular solution of~\eqref{eq:radial-ode} satisfying~\eqref{eq:radial-bc}, with center slope $\alpha=-\theta'(0)$. Under the parameter condition~\eqref{eq:range}, we have $\alpha>\sqrt h$.
\end{corollary}
\begin{proof}
The integrating factor $I$ and identity~\eqref{eq:IG} are valid for every regular strictly decreasing profile under $h>|A|$. Since $0<\nu(r)<1$ for $r>0$, condition~\eqref{eq:range} and the strict monotonicity in~\eqref{eq:R} give
\[
\kappa\le\cR(1)<\cR(\nu(r))
\qquad\text{for every }r>0.
\]
The second inequality is strict even when $4\kappa\sqrt h=3h-A$. Hence $(IG)'<0$ on $(0,\infty)$ by~\eqref{eq:IG}. Choosing a sufficiently large  $R$ so that $G(R)>0$ by~\eqref{eq:tail-signs}, we obtain
\[
2(\alpha-\sqrt h) =\lim_{r\downarrow0}I(r)G(r) >I(R)G(R)>0.
\]
Thus $\alpha>\sqrt h$.
\end{proof}

\section{Construction of the calibration}\label{sec:construction}
Assume $\kappa>0$ and $h>|A|$. Throughout this section, let $\theta$ be a regular solution of~\eqref{eq:radial-ode} satisfying~\eqref{eq:radial-bc}, with center slope $\alpha=-\theta'(0)>\sqrt h$.
Let $\M$ be the associated map~\eqref{eq:M}. For $x=r(\cos\varphi,\sin\varphi)$ with $r>0$, the image $\M(x)$ has spherical coordinates
\[
(t,\psi)=(\theta(r),\varphi).
\]
Let $r=R(t)$ be the inverse of the strictly decreasing function $t=\theta(r)$. At the target point
\[
\n(t,\psi)= \begin{pmatrix}
\sin t \cos \psi \\ \sin t \sin \psi \\ \cos t
\end{pmatrix}
, \qquad 0<t<\pi, \quad  \psi \in \R / 2\pi\Z,
\]
we use the positively oriented orthonormal frame
\[
\ein_t = \p_t \n=\begin{pmatrix}
\cos t \cos \psi \\ \cos t \sin \psi \\ -\sin t
\end{pmatrix}, \qquad
\ein_\psi=\frac{1}{\sin t} \p_{\psi} \n=\begin{pmatrix}
-\sin \psi \\ \cos \psi \\ 0
\end{pmatrix}.
\]
For the planar frame of Section~\ref{sec:intro}, we use the angular notation $\ein_r(\psi)=(\cos\psi,\sin\psi)$ and $\ein_\varphi(\psi)=R_{\pi/2}\ein_r(\psi)$.
For a target point $\n=\n(t,\psi)$ away from the poles, let $X(\n)=R(t)(\cos \psi, \sin \psi)$ be its preimage under $\M$. Define the reference Jacobian field by
\[
\fJ(\n):=D\M(X(\n)): \R^2 \longrightarrow T_{\n}\Sph.
\]
In this section, $a$ and $b$ denote the singular values $a(r)$ and $b(r)$ from Section~\ref{sec:barriers} evaluated at $r=R(t)$. In the planar and spherical frames above, the matrix of $\fJ$ is
\begin{equation}\label{eq:J}
\fJ(\n)=\operatorname{diag}(-a,b).
\end{equation}
In the same frames, the interfacial DMI coefficient $G(\n)$ and the adjoint projection differential $(\dif\pi(\n))^{\fT}$ are represented by
\[
G(\n)=\operatorname{diag}(1,\cos t), \qquad (\dif\pi(\n))^{\fT}=\operatorname{diag}(\cos t,1).
\]
For the pointwise construction, we use the energy density with the interfacial Lifshitz term
\[
f(\n,B)=\frac12|B|^2+\kappa G(\n):B+W(\n),
\qquad B:\R^2\to T_{\n}\Sph.
\]
The graph pullback formula~\eqref{eq:graph-pullback} suggests the candidate
\[
\ell(\n,B)=p(\n):B+\beta(\n)\det B.
\]
For the whole-plane comparison, we subtract $\kappa(\dif\pi(\n))^{\fT}:B$ from both $f$ and $\ell$. This replaces $G$ by $\widetilde G$ in the energy density as in~\eqref{eq:gauged-energy} and leaves the pointwise difference $f-\ell$ unchanged.

The weighted determinant gives the degree-dependent part of the calibration integral, and we seek to write the linear term as a divergence. To make the pointwise bound sharp at the reference gradient, we match the values and first $B$-derivatives of $f$ and $\ell$ at $B=\fJ(\n)$. The calculation leading to~\eqref{eq:criterion-data}, with $\cW=W$, $\fL=\kappa G$, and $p$ in place of $(\dif Y)^{\fT}$, gives the equations for $p$ and $\beta$.

Since $\det\fJ=-ab\ne0$ for $0<t<\pi$ and $|\fJ|^2=a^2+b^2$, the value identity determines
\begin{equation}\label{eq:beta}
\beta(t)=\frac{2W(t)-(a^2+b^2)}{2ab}.
\end{equation}
The singular-value bounds in Propositions~\ref{prop:lower} and~\ref{prop:upper} therefore give
\begin{equation}\label{eq:beta-gap}
1+\beta=\frac{2W-(a-b)^2}{2ab}>0,\qquad
1-\beta=\frac{(a+b)^2-2W}{2ab}>0.
\end{equation}
The derivative identity gives $p=\fJ+\kappa G-\beta\cof\fJ$. Since $\cof\fJ=\operatorname{diag}(b,-a)$ in the aligned frames, $p$ is diagonal in these frames, with entries
\begin{equation}\label{eq:p-components}
p_1=-a+\kappa-\beta b,\qquad
p_2=b+\kappa\cos t+\beta a.
\end{equation}
To construct a primitive for $p$, we seek $Y$ in the form $Y(\n)=y(t)\ein_r(\psi)$. Since $\ein_\psi=(\sin t)^{-1} \p_\psi \n$, the angular entry of $\dif Y$ is $y(t)/\sin t$. Matching this entry with $p_2(t)$ gives $y(t)=p_2(t)\sin t$, hence we introduce
\[
Y(\n)=p_2(t)n'.
\]
The remaining condition for $p=(\dif Y)^{\fT}$ is $(p_2\sin t)_t=p_1$. In the following lemma, we derive this identity from the radial equation and establish the pole extensions and decay estimates needed for the whole-plane comparison.
\begin{lemma}\label{lem:closedness}
 The fields $\fJ$, $\beta$ and $p$ extend continuously to the sphere. The map $Y$ extends to $C^1(\Sph;\R^2)$ and satisfies
\[
 (\dif Y)^{\fT}=p, \qquad \norm{\beta}_{L^\infty(\Sph)}<1.
\]
 After setting $\widetilde Y=Y-\kappa\pi$ and $\widetilde p=p-\kappa(\dif\pi)^{\fT}$, one has
\begin{equation}\label{eq:decay-data}
 |\fJ(\n)|+|\widetilde p(\n)|\le C|\n-\ein_3|,\qquad
 |\widetilde Y(\n)|\le C|\n-\ein_3|^2.
\end{equation}
\end{lemma}
\begin{proof}
Define $v(t)=-R(t)\theta'(R(t))>0$.  The inverse-function relation and the radial equation~\eqref{eq:radial-ode} give
\begin{equation}\label{eq:vt}
R_t = -\dfrac{R}{v}, \qquad v_t =-(r\theta')' \big|_{r=R(t)} R_t = \frac{1}{v}\bra*{\sin t \cos t - 2\kappa R \sin^2 t + R^2 W_t}.
\end{equation}
It follows from~\eqref{eq:beta} and~\eqref{eq:p-components} that
\[
p_1 = \kappa + \frac{\sin^2 t - v^2 -2R^2 W}{2Rv}, \qquad
p_2 \sin t = \kappa\sin t\cos t +RW +\frac{\sin^2 t-v^2}{2R}.
\]
Differentiating the second expression and using~\eqref{eq:vt} gives
\[
\begin{split}
(p_2 \sin t)_t &= \kappa(\cos^2 t-\sin^2 t) - \frac{RW}{v}+RW_t +\frac{\sin t \cos t - v v_t}{R} + \frac{\sin^2 t - v^2}{2Rv}\\
&=\kappa+\frac{\sin^2 t-v^2}{2Rv}-\frac{RW}{v}=p_1.
\end{split}
\]
The cancellation of the $W_t$ terms uses only the radial Euler--Lagrange equation. With the same exchange and DMI terms, the compatibility identity holds for any smooth axisymmetric potential wherever the radial profile is smoothly invertible.
Writing $Y=(p_2\sin t)\ein_r$ gives $\dif Y(\ein_t)=p_1\ein_r$ and $\dif Y(\ein_\psi)=p_2\ein_\varphi$, hence $(\dif Y)^{\fT}=p$ away from the poles.

At the south pole $-\ein_3$, $t \to \pi$ and $r=R(t) \to 0$. Set $q=\pi-t$. The center expansion gives
\[
R=\frac{q}{\alpha}+O(q^3), \qquad a=\alpha+O(q^2), \qquad b=\alpha+O(q^2).
\]
Consequently,
\[
\beta(\pi)=\frac{2h}{\alpha^2}-1 \in (-1,1), \qquad p_1(\pi)=\kappa-\frac{2h}{\alpha}=-p_2(\pi).
\]
In the horizontal chart $\xi=\pi(\n)=n'$ near the south pole, \eqref{eq:J} gives
\[
\dif \pi \circ \fJ=bI+(-a\cos t-b)\ein_r\otimes\ein_r.
\]
As $t \to \pi$, $b \to \alpha$ and $-a\cos t-b\to0$. The horizontal matrix converges uniformly in $\psi$ to $\alpha I$, so $\fJ$ extends with $\fJ(-\ein_3)\xi=(\alpha\xi,0)$. With the outward orientation, $\det\fJ(-\ein_3)=-\alpha^2$.

At the north pole $\ein_3$, $t \to 0$ and $r=R(t) \to \infty$. With $H=\sqrt{h+A}$, relation~\eqref{eq:tail_derivative} gives
\[
a=H\theta + \frac{\theta}{2r} + O\bra*{\frac{\theta}{r^2}}, \qquad b=\frac{\theta}{r}+O\bra*{\frac{\theta^3}{r}}.
\]
Thus
\[
\begin{split}
& W=\frac{H^2\theta^2}{2} + O(\theta^4), \qquad ab=\frac{H\theta^2}{r}+o\bra*{\frac{\theta^2}{r}},\\
& \frac{a^2}{2} = \frac{H^2\theta^2}{2}+\frac{H\theta^2}{2r}+o\bra*{\frac{\theta^2}{r}}, \qquad \frac{b^2}{2} = o\bra*{\frac{\theta^2}{r}}.
\end{split}
\]
The decay~\eqref{eq:tail_K} ensures that $O(\theta^4)=o(\theta^2/r)$ here. Thus $\beta(t)\to-1/2$, $p_1,p_2\to\kappa$, and $\fJ\to0$. We extend $\beta$ by $\beta(0)=-1/2$ and $\fJ$ by $\fJ(\ein_3)=0$. The two pole values of $\beta$ and~\eqref{eq:beta-gap} give $\norm{\beta}_{L^\infty(\Sph)}<1$ by compactness. The formula $p=\fJ+\kappa G-\beta\cof\fJ$ extends $p$ continuously to both poles.

Set $Y=0$ at both poles and use the horizontal charts $\n=(\xi,\pm\sqrt{1-|\xi|^2})$ with $\xi=\pi(\n)$. Write $c=\kappa$ at the north pole and $c=p_2(\pi)$ at the south pole. Since $p_2(t)\to c$, we have
\[
Y(\xi)=p_2(t)\xi=c\xi+o(|\xi|),
\]
so $D_\xi Y(0)=cI$.
Away from $\xi=0$, the differential of the inverse chart sends $(\ein_r,\ein_\varphi)$ to $(\ein_t/\cos t,\ein_\psi)$.
Hence
\[
D_\xi Y=p_2I+\bra*{\frac{p_1}{\cos t}-p_2}\ein_r\otimes\ein_r \longrightarrow cI
\]
at either pole, using $p_1(\pi)=-p_2(\pi)$ at the south pole. Thus $Y$ is $C^1$ on the sphere, and $(\dif Y)^{\fT}=p$ holds at the poles by continuity.

Finally, near the north pole the tail estimates give $|\fJ|=O(t)$, and~\eqref{eq:p-components} gives $p_1-\kappa=O(t)$ and $p_2-\kappa=O(t)$. Thus $\widetilde p=O(t)$ and $\widetilde Y=(p_2-\kappa)n'=O(t^2)$. Since $t\le C|\n-\ein_3|$ locally, compactness away from the north pole yields~\eqref{eq:decay-data} globally.
\end{proof}

The contact formula~\eqref{eq:beta} and the south-pole calculation apply to every regular decreasing profile with $\alpha>0$. At the south pole,
\[
\beta(\pi)=\frac{2h}{\alpha^2}-1,\qquad
1-\beta(\pi)=2\bra*{1-\frac{h}{\alpha^2}}.
\]
Thus the quadratic remainder degenerates at the south pole when $\alpha=\sqrt h$ and has negative directions there when $\alpha<\sqrt h$.
Together with the sufficiency proved above, this establishes
\[
\norm{\beta}_{L^\infty(\Sph)}<1
\quad\Longleftrightarrow\quad \alpha>\sqrt h,
\]
which proves the uniform positivity criterion in Remarks~\ref{rem:intro}(1).

If a regular strictly decreasing profile has center slope $\alpha=\sqrt h$, the radial equation still gives the compatibility identity $(p_2\sin t)_t=p_1$. 
The pole extensions and decay estimates above also apply. By Propositions~\ref{prop:lower} and~\ref{prop:upper}, we have $|\beta(t)|<1$ for $0<t<\pi$. 
Together with the pole values $\beta(0)=-1/2$ and $\beta(\pi)=1$, this gives $|\beta|\le 1$ on $\Sph$, with strict inequality away from the south pole.
Applying the divergence correction and integrating~\eqref{eq:pointwise} using~\eqref{eq:graph-integral} gives $E(\m)\ge E(\M)$ for all $\m\in\cA_{-1}$. If equality holds, the remainder vanishes almost everywhere.
On $\{\m=-\ein_3\}$ we have $D\m=0$ almost everywhere, so the remainder there equals $W(-\ein_3)=2h>0$. This set therefore has measure zero. 
The strict bound on $|\beta|$ elsewhere gives $D\m=\fJ(\m)$ almost everywhere, and the rigidity argument in the proof of Proposition~\ref{prop:calibration} shows that $\m$ agrees almost everywhere with a translate of $\M$.
\begin{remark}\label{rem:stress-closure}
Stress conservation and local invertibility ensure that the two-form obtained from the contact construction is closed. Let $\m$ be a smooth critical map for a $C^2$ autonomous energy density $f=f(\n,B)$, and let $\mathcal{U}$ be an open set compactly contained in the domain of $\m$ such that $\m|_{\mathcal{U}}$ is a diffeomorphism onto its image. On $\m(\mathcal{U})$, set $X=(\m|_{\mathcal{U}})^{-1}$ and $\fJ=D\m\circ X$. We match $\ell(\n,B)=p(\n):B+\beta(\n)\det B$ to $f$ in value and first $B$-derivative at $B=\fJ(\n)$ and obtain
\[
p=f_B-\beta \cof\fJ \quad\text{and}\quad \beta\det\fJ=f_B:\fJ-f.
\]

Stationarity under inner variations yields
\[
\sum_i\p_iS_{ij}=0,
\qquad S_{ij}=f_{B_i}(\m,D\m)\cdot\p_j\m-f(\m,D\m)\delta_{ij},
\]
where $B_i$ denotes the $i$th column of $B$. The tensor $S$ is the configurational stress measuring the energy response to spatial rearrangements of the magnetic texture. 
This is the static Noether balance associated with spatial translation invariance. 
The contact conditions and the cofactor identity give $p^{\fT}\fJ = S-\mathrm{tr}(S)I$ along $\m$. 
Define the candidate target one-forms by 
\[
\omega_i(\n)[v]=\scp{p(\n)e_i, v} \quad\text{for}\quad v\in T_{\n}\Sph, \, i=1,2.
\]
Their pullbacks are the one-forms associated with the rotated stress columns
\[
\m^*\omega_1=-S_{22}\dif x_1 + S_{12}\dif x_2, \qquad
\m^*\omega_2=S_{21}\dif x_1 - S_{11}\dif x_2.
\]
Stress conservation therefore gives $\dif(\m^*\omega_i)=0$. Local invertibility then implies
\[
\dif\omega_i=X^*\dif(\m^*\omega_i)=0.
\]
Consequently, the two-form
\[
\Omega=\omega_1\wedge\dif x_2-\omega_2\wedge\dif x_1+\beta\sigma
\]
is closed on $\R^2\times\m(\mathcal U)$, since $\beta\sigma$ has top degree on the target.
Its graph pullback has density $\ell$, so $\ell$ is locally a null Lagrangian. A calibration additionally requires the pointwise inequality $\ell\le f$.

For the radial reference map $\M$, radial stationarity gives $(rS_{rr})'=S_{\varphi\varphi}$. The contact conditions in the aligned frames imply
\[
p_2(\theta)\sin\theta=-rS_{rr},
\qquad
S_{\varphi\varphi}=ap_1(\theta), \qquad a=-\theta'.
\]
With $t=\theta(r)$, these identities yield $(p_2\sin t)_t=p_1$, the compatibility condition in Lemma~\ref{lem:closedness}.
The pole extensions and the strict gap require the additional estimates established there.
\end{remark}
\begin{proof}[Proof of Theorem~\ref{thm:main}]
Under condition~\eqref{eq:range}, Proposition~\ref{prop:shooting} gives a regular profile $\theta$ whose associated map is $\M\in\cA_{-1}$.
Writing $\alpha=-\theta'(0)$, Corollary~\ref{coro:range} gives $\alpha>\sqrt h$, so the calibration construction above applies.
Strict monotonicity and the endpoint values make $\M$ a bijection from $\R^2$ to $\Sph\setminus\{\ein_3\}$. Its differential is invertible away from the origin by~\eqref{eq:J}, and at the origin the center expansion gives $D(\pi\circ\M)(0)=\alpha I$ in the south-pole chart. Thus the inverse function theorem gives a smooth inverse on $\Sph\setminus\{\ein_3\}$.

For the whole-plane energy~\eqref{eq:gauged-energy}, we use
$\cW=W$, $\fL=\kappa\widetilde{G}$, and the data
$\fJ$, $\beta$, and $\widetilde{Y}$ from Lemma~\ref{lem:closedness}.
The identities
\[
W=\frac{1}{2}|\fJ|^2 - \beta\det\fJ, \qquad (\dif \widetilde{Y})^{\fT} = \fJ + \kappa \widetilde{G}-\beta \cof \fJ
\]
also hold at the poles by continuity. By~\eqref{eq:criterion-data}, they give the contact conditions for the corrected densities. Proposition~\ref{prop:calibration} therefore gives
\[
E(\m)-E(\M) \ge \frac{\gamma}{2} \int_{\R^2} |D\m-\fJ(\m)|^2 \dd x, \qquad \gamma=1-\norm{\beta}_{L^\infty(\Sph)}>0,
\]
and every global minimizer of $E$ over $\cA_{-1}$ equals $\M(\cdot-c)$ almost everywhere for some $c\in\R^2$. The smooth representative attains $-\ein_3$ only at $c$, so the translation is unique.

If there were another profile satisfying the conditions of Theorem~\ref{thm:main}, the center and tail arguments would apply to it as well. Under condition~\eqref{eq:range}, Corollary~\ref{coro:range} gives $\alpha>\sqrt{h}$ for its center slope, and Propositions~\ref{prop:lower} and~\ref{prop:upper} give both singular-value barriers. Hence the second profile also admits the calibration just constructed and is a global minimizer. Classification relative to $\M$ then makes it a translate of $\M$. Both centered maps attain $-\ein_3$ only at the origin, so that translation is zero. This proves uniqueness of the profile without assuming it in the shooting argument.
%
\end{proof}

\section{Linear stability}\label{sec:coercivity}
We use the calibration to identify the kernel of the Hessian at the centered minimizer $\M$. A far-field estimate and local compactness give the spectral gap in Corollary~\ref{coro:coercivity}(i). We then combine this gap with translation modulation to prove the local energy estimate in part~(ii).

We write
\[
\scp{\bs\eta, \nabla\M}_{L^2} = \bra*{\scp{\bs\eta, \p_1\M}_{L^2}, \, \scp{\bs\eta, \p_2\M}_{L^2}}.
\]

\begin{samepage}
\begin{lemma}\label{lem:kernel-gap}
Assume the hypotheses of Theorem~\ref{thm:main}, and let $\M$ be the centered minimizer. The Hessian $\cH_{\M}$ is nonnegative on $H^1(\R^2;T_{\M}\Sph)$ and has kernel
\[
\mathrm{span}\{\p_1\M,\p_2\M\}.
\]
Moreover, there exists $\lambda>0$ such that
\[
\cH_{\M}(\bs\eta)\ge\lambda\norm{\bs\eta-\bs\eta_0}_{H^1}^2
\qquad\text{for every }\bs\eta\in H^1(\R^2;T_{\M}\Sph),
\]
where $\bs\eta_0$ is the $L^2$-orthogonal projection of $\bs\eta$ onto the translation modes.
\end{lemma}
\end{samepage}
\begin{proof}
By Theorem~\ref{thm:main}, $\M$ minimizes $E$ in $\cA_{-1}$. The normalized variations of Section~\ref{sec:secondVariation} converge to $\M$ in the affine $H^1$ space as $t\to0$. By continuity of the degree, they remain in $\cA_{-1}$ for small $|t|$, so $\M$ is a critical point of $E$. By~\eqref{eq:M_regularity}, $\M-\ein_3\in H^2(\R^2;\R^3)$, and the energy-difference identity~\eqref{eq:Edifference} applies.

We use the exact defect identity~\eqref{eq:criterion-defect} from Proposition~\ref{prop:calibration}, with $C_{\m}=D\m-\fJ(\m)$. For $\bs\eta\in C_c^\infty(\R^2;T_{\M}\Sph)$, let $\m_t$ be the normalized variation from Section~\ref{sec:secondVariation}. On the compact support of $\bs\eta$, the images of $\M$ and $\m_t$ stay in a compact subset of $\Sph\setminus\{\ein_3\}$ for small $|t|$. There $\fJ=D\M\circ X$ is smooth, and we denote its differential as a matrix-valued map on the sphere by $\dif\fJ$. The normalization gives $(\m_t-\M)/t\to\bs\eta$ in $H^1$, and $C_{\m_t}=0$ outside the support. Hence
\[
C_{\m_t}=tT_{\bs\eta}+o(t)\quad\text{in }L^2,
\quad\text{where}\quad
T_{\bs\eta}=D\bs\eta-\dif\fJ(\M)[\bs\eta].
\]
Since $\m_t\cdot C_{\m_t}e_i=0$ and $C_{\M}=0$, differentiation at $t=0$ gives $\M\cdot T_{\bs\eta}e_i=0$.
Thus $T_{\bs\eta}:\R^2\to T_{\M}\Sph$, and its intrinsic determinant is well defined. By the scalar triple product formula, the strong $L^2$ convergence of $t^{-1}C_{\m_t}$ and the uniform convergence of $\m_t$ and $\beta(\m_t)$ give $t^{-2}\beta(\m_t)\det C_{\m_t}\to\beta(\M)\det T_{\bs\eta}$ in $L^1$.
Comparing the quadratic terms in~\eqref{eq:Edifference} and~\eqref{eq:criterion-defect} gives
\begin{equation}\label{eq:Heta}
\cH_{\M}(\bs\eta)
=\int_{\R^2}|T_{\bs\eta}|^2-2\beta(\M)\det T_{\bs\eta}\dd x
\ge\gamma\norm{T_{\bs\eta}}_{L^2}^2, \qquad
\gamma=1-\norm{\beta}_\infty>0.
\end{equation}
By~\eqref{eq:M_regularity}, the product rule and $H^1\hookrightarrow L^4$ make $P_{\M}$ bounded on $H^1$. Projecting smooth compactly supported approximations and using continuity gives $\cH_{\M}(\bs\eta) \ge 0$ for all $\bs\eta\in H^1(\R^2;T_{\M}\Sph)$.

Suppose $\cH_{\M}(\bs\eta)=0$. Applying~\eqref{eq:Heta} to a smooth tangent approximation and using local boundedness of $\dif\fJ(\M)$ yields $T_{\bs\eta}=0$ on each fixed ball. For such $\bs\eta$, set $w=\dif X(\M)[\bs\eta]\in H^1_{\rm{loc}}(\R^2;\R^2)$, where $X=\M^{-1}$. Since $\bs\eta=D\M w$ and $\fJ(\M(x))=D\M(x)$, differentiation gives
\[
T_{\bs\eta}=D\M Dw.
\]
For each $x\in\R^2$, $D\M(x):\R^2 \to T_{\M(x)}\Sph$ is a linear isomorphism. This gives $Dw=0$, so $w$ is constant and $\bs\eta$ is a linear combination of $\p_j\M$, $j=1,2$. Conversely, let $\M _t=\M (\cdot+t e_j)$. By translation invariance, $E(\M_t)=E(\M)$, and $(\M_t-\M)/t\to\p_j\M $ in $H^1$ since $\M -\ein_3\in H^2$. Dividing \eqref{eq:Edifference} by $t^2$ and passing to the limit gives $\cH_{\M}(\partial_j\M )=0$. Thus 
\begin{equation}\label{eq:kernel}
\cH_{\M}(\bs\eta,\bs\zeta)=0 \text{ for every } \bs\zeta\in H^1(\R^2;T_{\M}\Sph) 
\quad\text{if and only if}\quad
\bs\eta \in \mathrm{span}\set*{\p_1\M,\p_2\M}.
\end{equation}
For $\bs\eta\in H^1(\R^2;T_{\M}\Sph)$, the Hessian reads
\[
\cH_{\M}(\bs\eta) = \int_{\R^2} |D\bs\eta|^2 - 4\kappa \eta' \cdot \nabla \eta_3 + (h+A)|\eta'|^2 + h\eta_3^2 -\Lambda(\M)|\bs\eta|^2\dd x.
\]
Since $\M \to \ein_3$, we can take $R$ large enough that $M_3\ge\frac{1}{2}$ on $\R^2\setminus B_R$. There $\eta_3 = -M' \cdot \eta'/M_3$, and hence
\[
|\eta_3| \le C\delta_R|\bs\eta|, \quad |D\eta_3| \le C\delta_R(|D\bs\eta|+|\bs\eta|), \quad
\delta_R = \sup_{|x| \ge R}\bra*{|M'|+|D\M|+|\Lambda(\M)|} \to 0.
\]
Since $h+A>0$, the DMI term and the $\Lambda(\M)$ term can be absorbed into the positive terms on $\R^2\setminus B_R$ for large $R$. On $B_R$, the integrand is bounded below by $\frac{1}{2}|D\bs\eta|^2-C_R|\bs\eta|^2$ by Young’s inequality. Together these bounds give
\begin{equation}\label{eq:tangent-garding}
 \cH_{\M}(\bs\eta)  \ge c_0\norm{\bs\eta}_{H^1(\R^2)}^2-C_R\norm{\bs\eta}_{L^2(B_R)}^2.
\end{equation}
Suppose that $\cH_{\M}$ is not coercive in $H^1$ on the tangent fields
that are $L^2$-orthogonal to the translation modes.
Then there are $\bs\eta_k\in H^1(\R^2;T_{\M}\Sph)$
such that
\[
\scp{\bs\eta_k,\nabla\M}_{L^2}=0,\qquad
\norm{\bs\eta_k}_{H^1}=1,\qquad
\cH_{\M}(\bs\eta_k)\to0.
\]
After passing to a subsequence,
$\bs\eta_k\rightharpoonup\bs\eta$ in $H^1(\R^2)$
and $\bs\eta_k\to\bs\eta$ in $L^2(B_R)$.
The tangent constraint and the orthogonality conditions
pass to the weak limit.
By~\eqref{eq:tangent-garding},
\[
C_R\norm{\bs\eta_k}^2_{L^2(B_R)} \ge c_0 -\cH_{\M}(\bs\eta_k) \to c_0,
\]
so $\bs\eta\ne0$.
The nonnegative bounded quadratic form $\cH_{\M}$
is weakly lower semicontinuous on the closed tangent
subspace. Hence
\[
0\le\cH_{\M}(\bs\eta)
\le\liminf_{k\to\infty}\cH_{\M}(\bs\eta_k)=0.
\]
By~\eqref{eq:kernel}, $\bs\eta$ is a linear combination of the translation modes. Its orthogonality to those modes gives $\bs\eta=0$, a contradiction. This proves the $H^1$ coercivity on the $L^2$-orthogonal complement of the translation modes, with some $\lambda>0$.

For any $\bs\eta\in H^1(\R^2;T_{\M}\Sph)$, let $\bs\eta_0$ be its
$L^2$-orthogonal projection onto the translation modes. By~\eqref{eq:kernel}
and the coercivity bound for $\bs\eta-\bs\eta_0$,
\[
\cH_{\M}(\bs\eta)=\cH_{\M}(\bs\eta-\bs\eta_0)
\ge \lambda\norm{\bs\eta-\bs\eta_0}_{H^1}^2.
\]
This proves the claim.
\end{proof}

\begin{proof}[Proof of Corollary~\ref{coro:coercivity}]
Part~(i) follows from Lemma~\ref{lem:kernel-gap}.

\medskip
\noindent\emph{(ii) Local energetic stability.}
Assume $d_{H^1}(\m,\cO)<\eps$, where $\eps>0$ will be chosen sufficiently small. Choose $c_0\in\R^2$ such that $\norm{\m-\M_{c_0}}_{H^1}<\eps$, and set $\m_{c_0}=\m(\cdot+c_0)$.

Following~\cite[Lemma 10]{LiMelcher2018}, consider
\[
F:\R^2\times H^1(\R^2;\R^3)\to\R^2,
\qquad
F(a,v)=\scp{v(\cdot+a)-(\M-\ein_3),\nabla\M}_{L^2}.
\]
Its components are
\[
F_j(a,v)=\scp{v,\p_j\M(\cdot-a)}_{L^2} -\scp{\M-\ein_3,\p_j\M}_{L^2}.
\]
Since $\p_j\M\in H^1$, its translations depend $C^1$ on $a$ with values in $L^2$.
Thus $F$ is $C^1$. Moreover, $F(0,\M-\ein_3)=0$, and
\[
D_aF(0,\M-\ein_3)=\bigl(\scp{\p_i\M,\p_j\M}_{L^2}\bigr)_{i,j=1}^2
\]
is positive definite because the translation modes are linearly independent.
The implicit function theorem gives $a\in\R^2$ such that
\[
F(a,\m_{c_0}-\ein_3)=0, \qquad
|a|\le C\norm{\m_{c_0}-\M}_{H^1}.
\]
Set $c=c_0+a$ and $\bs\xi=\m(\cdot+c)-\M$. Then
\[
\scp{\bs\xi,\nabla\M}_{L^2}=0,
\qquad \norm{\bs\xi}_{H^1}
\le\norm{\m_{c_0}-\M}_{H^1}+\norm{\M(\cdot+a)-\M}_{H^1}
\le C\eps.
\]
The last estimate follows from $\M-\ein_3\in H^2$. All constants are independent of $c_0$, since the implicit function theorem is applied near the same point $(0,\M-\ein_3)$ after translation.

Moreover, translation invariance and~\eqref{eq:Edifference} yield
\[
E(\m)-E(\M)=\frac{1}{2}\cH_{\M}(\bs\xi).
\]
Decompose $\bs\xi=\bs\xi^\top+\bs\xi^\perp$, where
\[
\bs\xi^\top = P_{\M}\bs\xi \quad\text{and}\quad \bs\xi^\perp = (\M\cdot \bs\xi)\M=-\frac{1}{2}|\bs\xi|^2 \M.
\]
The embedding $H^1(\R^2)\hookrightarrow L^4(\R^2)$ and boundedness of $P_{\M}$ on $H^1$ give
\[
\norm{\bs\xi^\perp}_{L^2}=\frac12\norm{\bs\xi}_{L^4}^2 \le C\norm{\bs\xi}_{H^1}^2,
\qquad
\norm{\bs\xi^\top}_{H^1} \le C\norm{\bs\xi}_{H^1}.
\]
Since $\bs\xi^\perp$ is parallel to $\M$ and $\M\cdot\p_j\M=0$, we have
\[
\scp{\bs\xi^\top,\p_j\M}_{L^2}
=\scp{\bs\xi,\p_j\M}_{L^2}=0, \qquad j=1,2.
\]
Thus the $L^2$ projection of $\bs\xi^\top$ onto the translation modes vanishes. Combining Lemma~\ref{lem:kernel-gap} with~\eqref{eq:Edifference} gives
\[
E(\m)-E(\M)
\ge\frac12\bra*{\lambda\norm{\bs\xi^\top}_{H^1}^2+\cH_{\M}(\bs\xi^\perp)+2\cH_{\M}(\bs\xi^\top,\bs\xi^\perp)}.
\]
Put $f=-|\bs\xi|^2/2$, so that $\bs\xi^\perp=f\M$. Since $|\bs\xi|\le2$ and $\bs\xi\in H^1$, the field $f^2\M\in H^1$. 
Using it in~\eqref{eq:Hinfity}, we write the normal energy as the difference of two $\cH_\infty$ terms. Their DMI contributions cancel pointwise by the product rule
\[
2f M'\cdot\nabla(fM_3)
= M'\cdot\nabla(f^2M_3)+f^2 M'\cdot\nabla M_3.
\]
Expanding the exchange and potential terms then gives
\begin{align*}
\cH_{\M}(\bs\xi^\perp)
&=\cH_\infty(f\M)-\cH_\infty(\M-\ein_3,f^2\M)\\
&=\int_{\R^2}
\bra*{|\nabla f|^2+hM_3f^2}\dd x
=\int_{\R^2}
\bra*{|D\bs\xi^\perp|^2+(hM_3-|D\M|^2)f^2}\dd x\\
&\ge\norm{\bs\xi^\perp}_{H^1}^2
-C\norm{\bs\xi}_{H^1}^4.
\end{align*}
For the mixed exchange term, differentiating $\bs\xi^\top\cdot\M=0$ gives
\[
D\bs\xi^\top:D(f\M)
= fD\bs\xi^\top:D\M-\sum_i(\p_i f)\bs\xi^\top\cdot\p_i\M.
\]
Using $|Df|\le|\bs\xi||D\bs\xi|$, $|D\bs\xi^\top|\le C(|D\bs\xi|+|\bs\xi|)$, and boundedness of $\M$ and $D\M$, we bound the mixed exchange, DMI, and zero-order integrands by $C(|\bs\xi|^2|D\bs\xi|+|\bs\xi|^3)$. Hölder's inequality and the embeddings of $H^1$ into $L^3$ and $L^4$ yield
\[
\abs*{\cH_{\M}(\bs\xi^\top,\bs\xi^\perp)} \le C\bra*{\norm{D\bs\xi}_{L^2}\norm{\bs\xi}_{L^4}^2+\norm{\bs\xi}_{L^3}^3} \le C\norm{\bs\xi}_{H^1}^3.
\]
Combining the estimates with $\mu=\min\{\lambda,1\}$ and using $\norm{\bs\xi}_{H^1}^2\le 2(\norm{\bs\xi^\top}_{H^1}^2+\norm{\bs\xi^\perp}_{H^1}^2)$ gives
\[
E(\m)-E(\M)\ge\frac{\mu}{4}\norm{\bs\xi}_{H^1}^2
-C\bigl(\norm{\bs\xi}_{H^1}^3+\norm{\bs\xi}_{H^1}^4\bigr).
\]
Since $\norm{\bs\xi}_{H^1}\le C\eps$, we can take $\eps$ small enough to absorb the last two terms and obtain
\[
E(\m)-E(\M) \ge c_1\norm{\bs\xi}_{H^1}^2 \ge c_1d_{H^1}(\m,\mathcal O)^2.
\]
This proves the lower bound in~\eqref{eq:local-coercivity}.
For any $c\in\R^2$, apply~\eqref{eq:Edifference} to $\m(\cdot+c)$. Translation invariance and boundedness of the Hessian give
\[
E(\m)-E(\M)=\frac{1}{2}\cH_{\M}(\m(\cdot+c)-\M) \le C\norm{\m-\M_c}^2_{H^1}.
\]
Taking the infimum over $c\in\R^2$ proves the upper bound in~\eqref{eq:local-coercivity}.
\end{proof}

\bibliographystyle{abbrv}
\bibliography{references}

@article{BartonSinger2020,
  author  = {Barton-Singer, B. and Ross, C. and Schroers, B. J.},
  title   = {Magnetic skyrmions at critical coupling},
  journal = {Communications in Mathematical Physics},
  volume  = {375},
  pages   = {2259--2280},
  year    = {2020},
  doi     = {10.1007/s00220-019-03676-1},
  url     = {https://doi.org/10.1007/s00220-019-03676-1}
}

@article{BelavinPolyakov1975,
  author  = {Belavin, A. A. and Polyakov, A. M.},
  title   = {Metastable states of two-dimensional isotropic ferromagnets},
  journal = {JETP Letters},
  volume  = {22},
  pages   = {245--247},
  year    = {1975}
}

@article{BernandMantel2020,
  author  = {Bernand-Mantel, A. and Muratov, C. B. and Simon, T. M.},
  title   = {Unraveling the role of dipolar versus {Dzyaloshinskii--Moriya} interactions in stabilizing compact magnetic skyrmions},
  journal = {Physical Review B},
  volume  = {101},
  pages   = {045416},
  year    = {2020},
  doi     = {10.1103/PhysRevB.101.045416},
  url     = {https://doi.org/10.1103/PhysRevB.101.045416}
}

@article{BernandMantel2021,
  author  = {Bernand-Mantel, A. and Muratov, C. B. and Simon, T. M.},
  title   = {A quantitative description of skyrmions in ultrathin ferromagnetic films and rigidity of degree {$\pm1$} harmonic maps from {$\mathbb{R}^2$} to {$\mathbb{S}^2$}},
  journal = {Archive for Rational Mechanics and Analysis},
  volume  = {239},
  pages   = {219--299},
  year    = {2021},
  doi     = {10.1007/s00205-020-01575-7},
  url     = {https://doi.org/10.1007/s00205-020-01575-7}
}

@article{BrezisNirenberg1995,
  author  = {Brezis, H. and Nirenberg, L.},
  title   = {Degree theory and {BMO}; {Part I}: compact manifolds without boundaries},
  journal = {Selecta Mathematica. New Series},
  volume  = {1},
  pages   = {197--263},
  year    = {1995},
  doi     = {10.1007/BF01671566},
  url     = {https://doi.org/10.1007/BF01671566}
}

@misc{DiFrattaMuratovSimonSlastikov2026,
  author        = {Di Fratta, G. and Muratov, C. B. and Simon, T. M. and Slastikov, V. V.},
  title         = {Uniqueness and radial symmetry of magnetic skyrmions},
  year          = {2026},
  howpublished  = {arXiv:2611.02141v1},
  eprint        = {2611.02141},
  archivePrefix = {arXiv},
  url           = {https://arxiv.org/abs/2611.02141}
}

@article{DoeringMelcher2017,
  author  = {D{\"o}ring, L. and Melcher, C.},
  title   = {Compactness results for static and dynamic chiral skyrmions near the conformal limit},
  journal = {Calculus of Variations and Partial Differential Equations},
  volume  = {56},
  pages   = {60},
  year    = {2017},
  doi     = {10.1007/s00526-017-1172-2},
  url     = {https://doi.org/10.1007/s00526-017-1172-2}
}

@article{GustafsonWang2021,
  author  = {Gustafson, S. and Wang, L.},
  title   = {Co-rotational chiral magnetic skyrmions near harmonic maps},
  journal = {Journal of Functional Analysis},
  volume  = {280},
  pages   = {108867},
  year    = {2021},
  doi     = {10.1016/j.jfa.2020.108867},
  url     = {https://doi.org/10.1016/j.jfa.2020.108867}
}

@article{Hattori2024,
  author  = {Hattori, K.},
  title   = {Smooth maps minimizing the energy and the calibrated geometry},
  journal = {Journal of Geometric Analysis},
  volume  = {34},
  pages   = {9},
  year    = {2024},
  doi     = {10.1007/s12220-023-01452-1},
  url     = {https://doi.org/10.1007/s12220-023-01452-1}
}

@article{HoffmannEtAl2017,
  author  = {Hoffmann, M. and Zimmermann, B. and M{\"u}ller, G. P. and Sch{\"u}rhoff, D. and Kiselev, N. S. and Melcher, C. and Bl{\"u}gel, S.},
  title   = {Antiskyrmions stabilized at interfaces by anisotropic {Dzyaloshinskii--Moriya} interactions},
  journal = {Nature Communications},
  volume  = {8},
  pages   = {308},
  year    = {2017},
  doi     = {10.1038/s41467-017-00313-0},
  url     = {https://doi.org/10.1038/s41467-017-00313-0}
}

@misc{IbrahimMiuraRomanShimizu2025,
  author        = {Ibrahim, S. and Miura, T. and Rom{\'a}n, C. and Shimizu, I.},
  title         = {Phase transition thresholds and chiral magnetic fields of general degree},
  year          = {2025},
  howpublished  = {arXiv:2512.24598v1},
  eprint        = {2512.24598},
  archivePrefix = {arXiv},
  url           = {https://arxiv.org/abs/2512.24598}
}

@article{IbrahimShimizu2023,
  author  = {Ibrahim, S. and Shimizu, I.},
  title   = {Phase transition threshold and stability of magnetic skyrmions},
  journal = {Communications in Mathematical Physics},
  volume  = {402},
  pages   = {2627--2640},
  year    = {2023},
  doi     = {10.1007/s00220-023-04779-6},
  url     = {https://doi.org/10.1007/s00220-023-04779-6}
}

@misc{IbrahimShimizu2026,
  author        = {Ibrahim, S. and Shimizu, I.},
  title         = {Global perturbation of isolated equivariant chiral skyrmions from the {Bogomol'nyi} case},
  year          = {2025},
  howpublished  = {arXiv:2503.03209v2},
  eprint        = {2503.03209},
  archivePrefix = {arXiv},
  url           = {https://arxiv.org/abs/2503.03209v2}
}

@article{KomineasMelcherVenakides2023,
  author  = {Komineas, S. and Melcher, C. and Venakides, S.},
  title   = {Chiral magnetic skyrmions across length scales},
  journal = {New Journal of Physics},
  volume  = {25},
  pages   = {023013},
  year    = {2023},
  doi     = {10.1088/1367-2630/acb6e6},
  url     = {https://doi.org/10.1088/1367-2630/acb6e6}
}

@article{Leonov2016,
  author  = {Leonov, A. O. and Monchesky, T. L. and Romming, N. and Kubetzka, A. and Bogdanov, A. N. and Wiesendanger, R.},
  title   = {The properties of isolated chiral skyrmions in thin magnetic films},
  journal = {New Journal of Physics},
  volume  = {18},
  pages   = {065003},
  year    = {2016},
  doi     = {10.1088/1367-2630/18/6/065003},
  url     = {https://doi.org/10.1088/1367-2630/18/6/065003}
}

@article{LiMelcher2018,
  author  = {Li, X. and Melcher, C.},
  title   = {Stability of axisymmetric chiral skyrmions},
  journal = {Journal of Functional Analysis},
  volume  = {275},
  pages   = {2817--2844},
  year    = {2018},
  doi     = {10.1016/j.jfa.2018.01.019},
  url     = {https://doi.org/10.1016/j.jfa.2018.01.019}
}

@article{Melcher2014,
  author  = {Melcher, C.},
  title   = {Chiral skyrmions in the plane},
  journal = {Proceedings of the Royal Society A},
  volume  = {470},
  pages   = {20140394},
  year    = {2014},
  doi     = {10.1098/rspa.2014.0394},
  url     = {https://doi.org/10.1098/rspa.2014.0394}
}

@article{MonteilMuratovSimonSlastikov2023,
  author  = {Monteil, A. and Muratov, C. B. and Simon, T. M. and Slastikov, V. V.},
  title   = {Magnetic skyrmions under confinement},
  journal = {Communications in Mathematical Physics},
  volume  = {404},
  pages   = {1571--1605},
  year    = {2023},
  doi     = {10.1007/s00220-023-04864-w},
  url     = {https://doi.org/10.1007/s00220-023-04864-w}
}

@article{MuratovSimonSlastikov2025,
  author  = {Muratov, C. B. and Simon, T. M. and Slastikov, V. V.},
  title   = {Existence of higher degree minimizers in the magnetic skyrmion problem},
  journal = {Archive for Rational Mechanics and Analysis},
  volume  = {249},
  pages   = {61},
  year    = {2025},
  doi     = {10.1007/s00205-025-02131-x},
  url     = {https://doi.org/10.1007/s00205-025-02131-x}
}

@article{OlverSivaloganathan1988,
  author  = {Olver, P. J. and Sivaloganathan, J.},
  title   = {The structure of null {Lagrangians}},
  journal = {Nonlinearity},
  volume  = {1},
  pages   = {389--398},
  year    = {1988},
  doi     = {10.1088/0951-7715/1/2/005},
  url     = {https://doi.org/10.1088/0951-7715/1/2/005}
}
\end{document}